\documentclass[final]{siamltex}

\usepackage{amsmath}
\usepackage{amssymb}
\usepackage{url}
\usepackage{enumerate}
\usepackage{hyperref}
\usepackage{psfrag}
\usepackage{graphicx}
\usepackage{mathtools}
\usepackage{tikz}

\usetikzlibrary{decorations.pathreplacing}

\newcommand\nset{\mathbb{N}}
\newcommand\rset{\mathbb{R}}

\newtheorem{algorithm}{Algorithm}
\title{Local matching indicators for transport problems with concave
  costs%   To display SVN information
}

\author{Julie Delon\footnotemark[2], Julien Salomon\footnotemark[3]
  \and Andrei {Sobolevski\footnotemark[4]}$\ $\footnotemark[5]}

\begin{document}
\maketitle

\renewcommand{\thefootnote}{\fnsymbol{footnote}}
\footnotetext[2]{LTCI CNRS, T{\'e}l{\'e}com ParisTech,46 rue Barrault
F-75634 Paris cedex 13, France (julie.delon@enst.fr).} 
\footnotetext[3]{CEREMADE, UMR CNRS 7534, Universit\'e de
  Paris-Dauphine, Place du Mar\'echal De Lattre De Tassigny, F-75775
  Paris cedex 16, France (salomon@ceremade.dauphine.fr).} 
\footnotetext[4]{Institute for information transmission problems
  (Kharkevich Institute), 19 B. Karetny per., 127994 Moscow, Russia
  (sobolevski@iitp.ru).}
\footnotetext[5]{UMI 2615 CNRS ``Laboratoire J.-V.~Poncelet,'' 11
  B. Vlasievski per., 119002 Moscow, Russia.} 
\renewcommand{\thefootnote}{\arabic{footnote}} 

%\tableofcontents

%--------------------------------------------------------------------
%               SECTION I
%--------------------------------------------------------------------

\begin{abstract}
   In this paper, we introduce a class of local indicators that enable to
  compute efficiently optimal transport plans associated to arbitrary
  weighted distributions of $N$ demands and $M$ supplies in $\rset$ in the
case where the cost function is concave. Indeed, whereas this problem
can be solved linearly when the cost is a convex function of the
distance on the line (or more generally when the cost matrix between
points is a Monge matrix), to the best of our knowledge, no simple
solution has been proposed for concave costs, which are more realistic
in many applications, especially in economic situations. The problem
we consider may be unbalanced, in the sense that the weight of all the
supplies might be larger than the weight of all the demands. We show
how to use the local indicators hierarchically to solve the
transportation problem for concave costs on the line. 
%The computational  cost of the indicators is small and independent of $N$. A hierarchical use of them enables to obtain an efficient algorithm.
\end{abstract}

\section{Introduction}

The origins of optimal transportation go back to the late eighteen
century, when Monge~\cite{monge81memoire} published his
\textit{M\'emoire sur la th\'eorie des d\'eblais et des remblais}
(1781). The problem, which was rediscovered and further studied by
Kantorovich in the 1940's, can be described in the following
way. Given two probability distributions $\mu$ and  $\nu$ on X and $c$
a %nonnegative
 measurable cost function on $X \times X$, find a joint probability
 measure $\pi$ on $X \times X$ with marginals $\mu$ and  $\nu$ and
 which minimizes the transportation cost 
\begin{equation}
  \label{eq:transport-optimal}
  \int \int_{X \times X} c(x,y) d\pi(x,y).
\end{equation}
Probability measures $\pi$ with marginals $\mu$ and $\nu$ are called
% \textit{transportation plans} or
\textit{transport plans}. A transport plan that minimizes the cost~\eqref{eq:transport-optimal} is said to be \textit{optimal}. 

When the measures $\mu$ and $\nu$ are discrete (linear combinations of
Dirac masses), the problem can be recast as finite linear 
programming. For $N\geq 1$, consider two discrete distributions of mass, or
\emph{histograms}, 
given on $\rset^N$: $\{(p_i, s_i)\}$, which represents ``supplies'' at locations $p_i$ with weights $s_i$
and $\{(q_j, d_j)\}$, which represents ``demands'' at locations $q_j$ with weights $d_j$
(notation
from~\cite{aggarwal1992efficient}) %
and assume that all values of $s_i$ and~$d_j$ are positive reals with  $S := \sum_i s_i$ and $D:= \sum_j d_j$. The problem consists in minimizing the transport cost
\begin{equation}
 \label{eq:transport_discret}
  \sum_{i,j} c(p_i, q_j) \gamma_{ij},
\end{equation}
where $\gamma_{ij}$ is the amount of mass going
from $p_i$ to~$q_j$, subject to the conditions
\begin{equation}
  \label{eq:transport-conditions}
  \gamma_{ij} \ge 0,\quad
  \sum_j \gamma_{ij} \le s_i,\quad
  \sum_i \gamma_{ij} \le d_j, \quad
\sum_{i,j} \gamma_{ij} =\min(S, D). 
\end{equation}
The matrix of values $\gamma = \{\gamma_{ij}\}$ is still called \emph{transport plan}.
 When $S=D$, the problem is said to be \textit{balanced} and is only a
 reformulation of~\eqref{eq:transport-optimal} for discrete
 measures. When~$S \neq D$, the problem
 is said to be \textit{unbalanced}. The cases $S<D$ and $S>D$ can be
 treated in the same way. This paper deals with balanced problems and
 unbalanced problems of the form $S>D$.

In the \textit{unitary case}, {\em i.e.} when all the masses $s_i$ and
$d_j$ are equal to a single value $v$, it turns out that if $\gamma$
is optimal, for all $i,j$, $\gamma_{ij} \in \{0,v\}$ and for all $ j$
there exists only one $i$ such that $\gamma_{ij} = v$ (each demand
receives all the mass from one supply). In the balanced case, the
matrix $\gamma$ is thus a permutation matrix up to the factor~$v$. In the unbalanced case,
the permutation matrix is padded with some zero rows. As a
consequence, the balanced
case boils down to an \textit{assignment problem}, known as the
\textit{linear sum assignment problem}. Such problems have been
thoroughly studied by the combinatorial optimization 
community~\cite{burkard2008assignment}. 
%and can be solved with classical linear programming methods, such as the Hungarian method or more involved algorithms (see~\cite{burkard2008assignment} for details).

Optimal transportation problems appear in many fields, such as economy
or physics for instance, see e.g.~\cite{brenier,cullen,lachapelle}. In
economic examples optimal transport is often related to the field of
logistic where
supplies are furnished by producers at specific places $p_i$ and in
specific quantities $d_i$, while demands corresponds to consumers
locations and needs. Depending 
on the application, various cost functions $c$ can be used. For
instance, concave functions of the distance appear as more realistic
cost functions 
in many economic situations. Indeed, as underlined by
McCann~\cite{mccann1999est}, a concave cost  
``\textit{translates into an economy of  scale for longer trips and may encourage cross-hauling.}''

\begin{figure}
\begin{center}
\begin{tikzpicture}
\draw (-2,0) --   (2,0) ;

\draw ( 1.2,0) arc  (0:180:1.2) ;
\draw[blue] ( 1.2,0) node {$\times$} ;
\draw[red] (-1.2,0) node {$\bullet$} ;

\draw ( 1,0) arc  (0:180:1) ;
\draw[blue] ( 1,0) node {$\times$} ;
\draw[red] (-1,0) node {$\bullet$} ;

\draw ( .8,0) arc  (0:180:.8) ;
\draw[blue] ( .8,0) node {$\times$} ;
\draw[red] (-.8,0) node {$\bullet$} ;

\draw ( .6,0) arc  (0:180:.6) ;
\draw[blue] ( .6,0) node {$\times$} ;
\draw[red] (-.6,0) node {$\bullet$} ;
\end{tikzpicture}
\hspace{.4cm}
\begin{tikzpicture}
\draw (-2,0) --   (2,0) ;

\draw ( 1.2,0) arc  (0:180:.9) ;
\draw[blue] ( 1.2,0) node {$\times$} ;
\draw[red] (-.6,0) node {$\bullet$} ;

\draw ( 1,0) arc  (0:180:.9) ;
\draw[blue] ( 1,0) node {$\times$} ;
\draw[red] (-.8,0) node {$\bullet$} ;

\draw ( .8,0) arc  (0:180:.9) ;
\draw[blue] ( .8,0) node {$\times$} ;
\draw[red] (-1,0) node {$\bullet$} ;

\draw ( .6,0) arc  (0:180:.9) ;
\draw[blue] ( .6,0) node {$\times$} ;
\draw[red] (-1.2,0) node {$\bullet$} ;
\end{tikzpicture}
\caption{On the left: optimal plan associated to a concave cost. On the right: optimal plan associated to a convex cost. Supplies are represented by red points and demands by blue crosses. }\label{excon}
\end{center}
\end{figure}
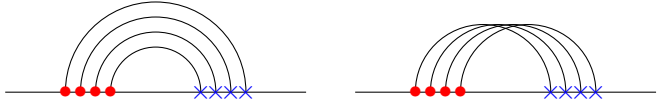

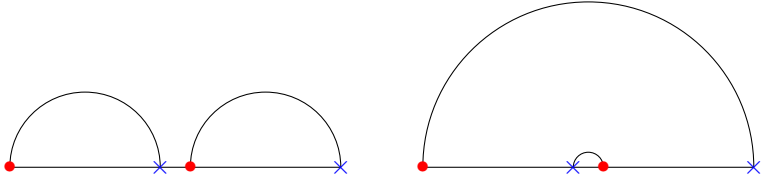
\begin{figure}
\begin{center}
\begin{tikzpicture}
\draw (-2,0) --   (2.4,0) ;

\draw ( 0,0) arc  (0:180:1) ;
\draw[blue] ( 0,0) node {$\times$} ;
\draw[red] (-2,0) node {$\bullet$} ;

\draw ( 2.4,0) arc  (0:180:1) ;
\draw[blue] ( 2.4,0) node {$\times$} ;
\draw[red] (0.4,0) node {$\bullet$} ;

\end{tikzpicture}
\hspace{.4cm}
\begin{tikzpicture}
\draw (-2,0) --   (2.4,0) ;

\draw ( 2.4,0) arc  (0:180:2.2) ;
\draw[blue] ( 2.4,0) node {$\times$} ;
\draw[red] (-2,0) node {$\bullet$} ;

\draw ( 0.4,0) arc  (0:180:.2) ;
\draw[blue] ( 0.0,0) node {$\times$} ;
\draw[red] (0.4,0) node {$\bullet$} ;

\end{tikzpicture}
\caption{On the left: solutions associated to the concave cost $c(x,y)
  = |x-y|^{0.9}$, and on the right to the cost $c(x,y) =
  |x-y|^{0.5}$. Supplies are represented by red points and demands by blue crosses. }\label{fig:contre_ex}
\end{center}
\end{figure}

During the last decades, many authors have taken interest in the study
of existence, uniqueness and properties of optimal
plans~\cite{Ambrosio03lecturenotes,mccann3,mccann2}, with a specific
interest for convex costs, {\em i.e.} costs $c$ that can be written as convex functions of the distance on the line. Detailed descriptions of these results can
be found in the books~\cite{villani2003tot,villani2008optimal}. One
case of particular interest is the one-dimensional case, which, when
$c$ is a convex function of the distance on the line, has been
completely understood~\cite{rachev1985monge} both for continuous and discrete settings. Indeed, this problem has
an explicit solution that does not depend on $c$ (provided that it is
convex) and consists in a monotone rearrangement (see Chapter 2.2
of~\cite{villani2003tot}).
In the unitary case,  this property can also be seen as a consequence of another interesting result, true for any dimension $N$, which says that the \textit{linear sum assignement problem} is solved by the identical permutation, provided that the cost matrix $(c(p_i,q_j))_{i,j}$ is a Monge matrix~\footnote{A matrix $C$ is said to be a Monge matrix if it satisfies $c_{ij} + c_{kl} \leq  c_{il} +c_{kj}$ when $ i <k$ and $ j <l$.}~\cite{burkard2008assignment}. 
Several approaches have been proposed to generalize the convex one-dimensional result to the case of the circle, where the starting point for the
monotone rearrangement is not known, and its choice and hence the optimal plan itself, unlike in the case of an interval, do depend on
the cost function $c$. Most of these approaches concern either the unitary case~\cite{karp1975two,werman_1985,werman1986bgm,cabrelli1995kantorovich,cabrelli1998linear,texture_metric} or the more general discrete case~\ref{eq:transport_discret}~\cite{pele2008eccv,ref:rabin08icprcemd,ref:rabin09siam,Pele_iccv2009, ref:rabin10jmiv}. Recently an efficient method has been introduced to tackle this issue in a continuous setting~\cite{dss10fast}. Unfortunately, these results on the line and the circle do not extend to
non-convex costs, in particular to concave costs (see
Figure~\ref{excon} for an example). Although it is of broad  interest
for many applications, few works treat this
case (see however the important paper~\cite{mccann1999est}) and computing solutions is far from obvious
in general. Indeed, contrary to the convex case on the line, optimal plans
strongly depend on the choice of the function $c$. Consider the case
of two unitary supplies at positions $p_1 = 0$ and $p_2 = 1.2$ and two
unitary demands at positions $q_1 = 1$ and $q_2=2.2$ on the line, as
drawn on Figure~\ref{fig:contre_ex}. If the cost function is $c(x,y) =
|x-y|^{0.9}$, the left solution will be optimal, whereas the other one
will be chosen for $c(x,y) = |x-y|^{0.5}$. For a convex cost, the left
solution would always be chosen. 

% The existence and properties of solutions for the transportation problem, continuous and discrete cases, McCann, Villani:
% \begin{itemize}
% \item explicit solution in 1D if $c$ is convex, rearrangement monotone, case of the circle~\cite{dss10fast}
% \item ND : no analytic solution
% \item non convex cost : nothing, even in 1D. 
% \item study of the concave case, 1D~\cite{mccann1999est}, \cite{aggarwal1992efficient} (linear case)
% \item unbalanced case : no explicit solution either in general.
% \end{itemize}

In practice, when no analytic solution is given ({\em i.e.} most of
the time), finding optimal plans can be a tedious task. As underlined
before, in a discrete setting, the problem can be written as a linear
programming problem, and optimal plans can be constructed numerically
by using for instance the simplex method or specialized methods such
as the auction algorithms~\cite{bertsekas1992auction} and various
algorithms for the assignment problem
(see~\cite{burkard2008assignment} for details). However these methods
do not take into account essential geometric features of the problem,
such as the fact that it is one-dimensional or that the cost function
is concave.

The goal of this paper is to introduce a class of functions that 
%can be used to 
%compute efficiently optimal transport plan in the
%one-dimensional case
reveals the local structure of optimal transport  plans
% in the 
%concave case
%The goal of this paper is to study the structural properties of
%optimal plans 
in the one-dimensional case, 
%(either on the line, or on the circle)
when the cost $c$ is a concave function of
the distance.
% The study will be first limited to the case of unitary 
% masses, either in the balanced or unbalanced case.
As a by-product, we
build an algorithm that permits to obtain optimal transport plans in the unitary case 
in less than $O(N^2)$ operations in both balanced and unbalanced cases, where $N$ is the
number of points under consideration. Once generalized to the non unitary case, the complexity of this algorithm becomes $O(N^3)$ in the worst case
but turns out to be smaller for ``typical'' problem instances.
However, let us insist that
our aim is not to compete with recent linear assignment algorithms,
which may be more interesting in practice, at least for balanced
problems, but rather to achieve a more complete understanding
of the internal structure of the assignment problem for concave costs on the line.

Observe that our algorithm complements the method suggested
by McCann~\cite{mccann1999est}, although the approach we follow here is
closer to the purely combinatorial approach
of~\cite{aggarwal1992efficient}. The results of this last work, in
which the cost $c(x, y) = |x - y|$ was considered, are
extended here to the general
framework of strictly concave cost functions. (Note that the very
special case considered in~\cite{aggarwal1992efficient} may be also
regarded as convex, which allows to apply the sorting algorithm on the
line or results of~\cite{dss10fast} on the circle.)

The paper is organized as follows: 
Section~\ref{sec:pb} is devoted to the presentation of the optimal transport problem. 
In Section~\ref{sec:chains}, we focus on transport problems on ``chains'' which are particular cases where
demands and supplies are alternated. In this framework, we  
present the main result of the paper, which states that consecutive
matching points in the optimal plan can be found thanks to local
indicators, independently of other points on the
line.  Thanks to the low number of
evaluations of the cost function required to apply the indicators, we
derive from this result a rather efficient algorithm to compute optimal
transport plans.
We then consider more general frameworks, namely general unitary
cases in Section~\ref{sec:unitary} and real-valued masses situations in
Section~\ref{sec:nonunitary}.
In Section~\ref{sec:prac}, we conclude with remarks on the
implementation of our algorithm and show that its complexity scales as $O(N^2)$ in the worst case.
Some technical proofs about this last result are given
in Appendix.

%--------------------------------------------------------------------
%               SECTION II
%--------------------------------------------------------------------

\section{The optimal transport problem}\label{sec:pb}
 This paper deals with the problem of finding an optimal transport
 plan in the case where the problem contains possibly more supplies
 than demands and the transport cost is strictly concave: the 
 larger the distance to cover is, the less the transport costs per
 unit distance, while the marginal cost (the derivative of the cost
 function) decreases monotonicaly. 

Consider two integers $M$,~$N$ and two sets of points $P = \{p_i \colon
i = 1, \dots, M\}$ and  $Q = \{q_i \colon i = 1, \dots, N\}$
 in $\rset$ that represent respectively the supply and
 demand locations.
 Let $s_i > 0$ be the capacity of $i$th supply and $d_j > 0$ the capacity
 of $j$th demand. We suppose that $S := \sum_i s_i \ge D := \sum_j
 d_j$, {\em i.e.} that the problem may be unbalanced.

 We deal with minimizing the cost 
\begin{equation}\label{eq:4}
  C(\gamma) = \sum_{i,j} c(p_i, q_j) \gamma_{ij},
\end{equation}
where $c(p_i, q_j)\in\rset^+ $ is the cost resulting
from transport of a unit mass between $p_i$ and $q_j$. The
quantity  $\gamma_{ij}$ is the amount of mass going from $p_i$ to
$q_j$, subject for all $i$,~$j$ to the conditions 
\begin{equation}
  \label{eq:5}
  \gamma_{ij} \ge 0,\quad
  \sum_j \gamma_{ij} \le s_i,\quad
  \sum_i \gamma_{ij} = d_j
\end{equation}
(observe that since $D\le S$, these conditions are equivalent to~\eqref{eq:transport-conditions}).
We call the case $S = D$ \textit{balanced} and the
case $S > D$ \textit{unbalanced}. Observe that
in the latter case the total supply is larger that the 
total demand, and therefore some of the supplies may remain
\emph{underused} ($\sum_j \gamma_{ij} < s_i$). 

As mentioned in Introduction, an optimal transport problem associated to equal masses,
i.e. $\forall (i,j)\in P\times Q$, $s_i=d_j=v$, reduces
actually to an assignment problem, where masses cannot be cut. Indeed, it is well known
(see Section 2.2 of~\cite{burkard2008assignment} for a proof) that if $\gamma$ minimizes the 
cost~\eqref{eq:4} under conditions~\eqref{eq:5}, then without loss
of generality one can assume that $\gamma_{ij} \in
\{0,v\}$ for all $i$,~$j$, so that the problem can be reformulated as finding
the minimum of the quantity
\begin{equation}\label{eq:1}
C(\sigma) = \sum_{1\le j\le N} c(p_{\sigma^{-1}(j)}, q_j),
\end{equation}
over all partial maps $\sigma\colon \{1, \dots, M\} \to \{1, \dots, N\}$
whose inverse~$\sigma^{-1}$ is injective and defined for all~$1\le
j\le N$: namely $j = \sigma(i)$ and $i = \sigma^{-1}(j)$ iff
$\gamma_{ij} = 1$. This setting is the one of
Sections~\ref{sec:chains} and~\ref{sec:unitary}.

We focus on the case where the function $c$ involves a strictly concave
function as stated in the next definition.
\begin{definition}\label{def:cost}
  The \emph{cost function} $c$ in~\eqref{eq:1} is said to be
  \emph{concave} if it is defined by
$c(p, q) = g(|p - q|)$ with $p, q \in \rset$, where
$g\colon \rset^+\rightarrow\rset \cup \{-\infty\}$ is a strictly concave non-decreasing 
  function such that $g(0) := \lim_{x\to 0} g(x) \ge -\infty$. 
\end{definition}

Note that strict concavity of~$g$ implies its \emph{strict}
monotonicity. 
Some examples of such costs are given by $g(x) =
\log x$ with $g(0) = -\infty$ and $g(x) = \sqrt x$ with
$g(0) = 0$.
If $g(0) > -\infty$, we assume without loss of generality that
$g(0) = 0$ (this changes the value of~(\ref{eq:4}) by an amount
$D\, g(0)$ independent of the transport plan).

In what follows, we denote by $\gamma^\star$ a given
optimal transport plan between $P$ and $Q$: $C(\gamma^\star)\leq C(\gamma)$
for all $\gamma$ satisfying~\eqref{eq:5}.
Observe that if two points $p_i$ and $q_j$ have the same position,
then there exists an optimal transport plan $\gamma^\star$ between $P$ and
$Q$ such that $\gamma^\star_{ij} = \min\{s_i, d_j\}$, {\em i.e.} that all
mass shared by the two marginal measures stays in place~\cite{villani2003tot}.
Indeed, suppose that a supply $p$
and a demand $q$  located at the same point are not matched together but
to some other demand and supply  $p'$ and $q'$ located at distances $x$ and $y$ respectively. 
Irrespective of whether $g(0)=0$ or
$g(0)=-\infty$, as soon as $g$ is strictly concave, one has
$$g(0) + g(x + y) < g(x) + g(y)$$
for all $x$,~$y$, 
which implies that matching $p$ and $q$ is cheaper.
Therefore a common point of $P$ and~$Q$ with unequal values $s_i$
and~$d_j$ may be replaced with a single supply of capacity~$s_i -
d_j$, if this quantity is positive, or with a single demand of
capacity~$d_j - s_i$.
In the following, we will therefore assume that common points do not
exist, {\em i.e.} that the sets $P$ and~$Q$ are disjoint.

Another significant feature of concave costs is that 
trajectories of mass elements under an optimal transport plan 
do not cross each other, as described by the following lemma.
% in the following sense.

\begin{lemma}[``Non-crossing rule'']
Consider two pairs of points $(p,q)$ and 
$(p',q')$ such that 
\begin{equation}
  \label{eq:cost-inequality}
c (p,q) + c(p',q') \le c (p',q) + c(p,q').
  \end{equation}
%It is easy to prove that the  
Then, the open intervals 
\begin{displaymath}
  I = (\min(p,q), \max(p,q)),\quad
  I' = (\min(p',q'), \max(p',q'))
\end{displaymath}
are \emph{nested}, in the sense that the following alternative holds:
 \begin{enumerate}
 \item  either $I\cap I'$ is empty,
\item or one of these intervals is a subset of the other.
\end{enumerate}
\label{lemma:non-crossing-rule}  
\end{lemma}

This result directly follows from the concavity of the cost
function and is often referred to as the ``non-crossing rule''
\cite{aggarwal1992efficient,mccann1999est}. The proof is based on the
same ideas as used in \cite{mccann1999est}. % (see also the proof of Lemma~\ref{lm:non-crossing} below).
  Essentially, the case $p < q' < q < p'$ and the similar
case with $p$'s and~$q$'s interchanged are ruled out in view
of~\eqref{eq:cost-inequality} by monotonicity of~$g$, whereas the case
$p < p' < q < q'$ and the symmetrical one are ruled out by
the strict concavity of~$g$.

In the unbalanced case, some supplies may lie outside all nested
segments.
\begin{definition}
  A point $r\in P\cup Q$ is said to be \emph{exposed} in the
  transport plan~$\gamma$ if $r\notin (\min(p_i,q_j),\max(p_i,q_j))$
  whenever $\gamma_{ij} > 0$.
\end{definition}

A sufficient condition for a supply to be exposed is given in the next statement.
\begin{lemma}
  In the unbalanced case all underused supplies are exposed in an
  optimal transport plan.
\label{lemma:unma-isol}  
\end{lemma}

Indeed, should an underused supply~$p_i$ belong to the interval
between $p_{i_0}$ and~$q_{j_0}$ such that $\gamma_{i_0j_0} > 0$, the
amount of mass equal to $\min\{\gamma_{i_0j_0}, s_i - \sum_j
\gamma_{ij}\}$ could be remapped to go to~$q_{j_0}$ from~$p_i$ rather
than~$p_{i_0}$, thus reducing the total cost of transport because of
the strict monotonicity of the function~$g$.

%--------------------------------------------------------------------
%               SECTION III
%--------------------------------------------------------------------

\section{Transport plans on chains}\label{sec:chains}
In this section, we focus on the particular case of {\it chains},
that are situations where all the masses are equal and alternated on the
line, i.e.  where $P$ and $Q$ satisfy $M=N$ ({\it balanced
  case}) and 
\begin{equation}\label{balanced}
 p_1 < q_1 < \dots < p_i < q_i < p_{i + 1} < q_{i
+ 1} < \dots < p_N < q_N,
\end{equation}
 or $M=N+1$ ({\it unbalanced case}) and 
\begin{equation}\label{eq:alter}
p_1 < q_1 < \dots < p_i < q_i < p_{i + 1} < q_{i
+ 1} < \dots < p_{N} < q_{N} < p_{N+1}.
\end{equation}
In these cases the set $P\cup Q$ is called {\it balanced chain} and {\it
  unbalanced chain} respectively. 
Sections coming after this one will aim at extenting our results to more
general cases.

Recall that in such a framework, optimal transport problems
are actually assignment problems, where masses cannot be cut: Optimal
transports plan are then described by permutations.

\subsection{Main result}
Thanks to the non-crossing rule, one knows that in any optimal
transport plan there exists at least two consecutive
points $(p_i,q_i)$ or $(q_i,p_{i+1})$ that are matched. Starting from this remark, we take advantage of the
structure of a chain to introduce
a class of indicators that enable to detect a priori such pairs of points.

\begin{definition}[Local Matching Indicators of order $k$]\label{def:LMI}
Given $0<k\leq N-1$, consider $2k+2$ consecutive points in a chain.
If the first point is a supply $p_i$, define 
$$
I^p_k(i)=c(p_i    ,q_{i+k})+\sum_{j=0}^{k-1}c(p_{i+j+1},q_{i+j})-\sum_{j=0}^{k}c(p_{i+j},q_{i+j}),
$$
else denote the first point $q_i$ and define
$$
I^q_k(i)=c(p_{i+k+1},q_{i})+\sum_{j=1}^{k}c(p_{i+j},q_{i+j})-\sum_{j=0}^{k}c(p_{i+j+1}  ,q_{i+j}).
$$
\end{definition}

\begin{figure}[h]
\begin{center}
\begin{tikzpicture}
\draw (0,0) --   (.7*7,0) ;

\draw[red] (.7*1,0) node {$\bullet$} ;
\draw[red] (.7*3,0) node {$\bullet$} ;
\draw[red] (.7*5,0) node {$\bullet$} ;

\draw[blue] (.7*2,0) node {$\times$} ;
\draw[blue] (.7*4,0) node {$\times$} ;
\draw[blue] (.7*6,0) node {$\times$} ;

\draw (.7*3,0) arc  (0:180:.7* .5) ;
\draw (.7*5,0) arc  (0:180:.7* .5) ;
\draw (.7*6,0) arc  (0:180:.7*2.5) ;
\end{tikzpicture}
\begin{tikzpicture}
\draw (0.5,1) node {\Huge $-$} ;
\end{tikzpicture}
\begin{tikzpicture}
\draw (0,0) --   (.7*7,0) ;

\draw[red] (.7*1,0) node {$\bullet$} ;
\draw[red] (.7*3,0) node {$\bullet$} ;
\draw[red] (.7*5,0) node {$\bullet$} ;

\draw[blue] (.7*2,0) node {$\times$} ;
\draw[blue] (.7*4,0) node {$\times$} ;
\draw[blue] (.7*6,0) node {$\times$} ;

\draw (.7*2,0) arc  (0:180:.7* .5) ;
\draw (.7*4,0) arc  (0:180:.7* .5) ;
\draw (.7*6,0) arc  (0:180:.7* .5) ;
\end{tikzpicture}
\caption{Schematic representation of an indicator of order 2.}\label{indic}
\end{center}
\end{figure}
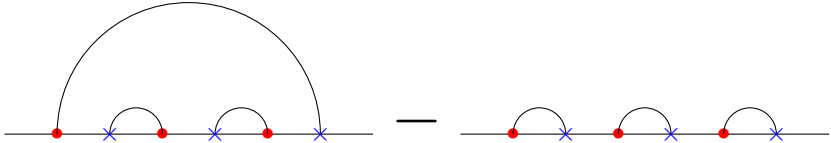

This definition is schematically depicted in Figure~\ref{indic} in the
case $k=2$.

Note that in the first alternative of this definition, we have
necessarily $1\leq k\leq N-1$,  $1\leq i\leq N-k$. 
In the second alternative, we have necessarily  $1\leq k\leq N-2$ and
$1\leq i\leq N-k-1$ in the balanced case and $1\leq k \leq N-1$ and
$1\leq i\leq N-k$ in the unbalanced case. 
The interest of these functions lies in the next result.

\begin{theorem}[Negative Local Matching Indicators of order $k$]\label{ThLMI}
%Assume that we are in the balanced case. 
Let $k_0\in \nset$ with $1\leq k_0\leq N-1$ and $i_0\in \nset$, %(resp. $i_0'\in \mathbf{N}$)  
such that $1\leq i_0\leq N-k_0$. In the unbalanced case, suppose in addition that
$g$ is strictly monotone.
%(resp. $1\leq i_0'\leq N-k_0-1$).

Assume that 
\begin{enumerate}
\item  $I^p_k(i) \geq 0$ for $k=1,\dots,k_0-1$, 
$i_0\leq i \leq i_0+k_0-k$,\label{hyp1} 
%was $1\leq i \leq N-k$,\label{hyp1} 
\item  $I^q_k(i')\geq 0$ for $k=1,\dots,k_0-1$,
$i_0\leq i'\leq i_0+k_0-k-1$, \label{hyp2}    
%was $1\leq i'\leq N-k-1$, \label{hyp2}   
(resp. $i_0\leq i' \leq i_0+k_0-k$ in the unbalanced case)
\item  $I^p_{k_0}(i_0) < 0$. %(resp. $I^q_{k_0}(i_0')< 0$). 
\label{hyp3}
\end{enumerate}

Then any permutation $\sigma$ associated to an optimal transport
plan satisfies
$\sigma(i)=i-1$ for $i=i_0+1,\dots,i_0+k_0$.

If the third condition is replaced by $I^q_{k_0}(i_0)< 0$ (with the same bounds on $k_0$ and $i_0$ in the unbalanced case, and with $1\leq k_0\leq N-2$ and $1\leq i_0\leq N-k_0-1$ in the balanced case), then  any permutation $\sigma$ associated to an optimal transport
plan satisfies $\sigma(i)=i$ for $i=i_0+1,\dots,i_0+k_0$.
\end{theorem}

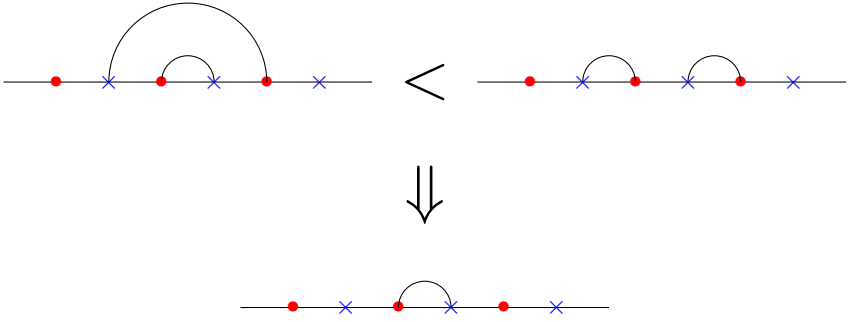
\begin{figure}
\begin{center}
\begin{tikzpicture}
\draw (0,0) --   (.7*7,0) ;

\draw[red] (.7*1,0) node {$\bullet$} ;
\draw[red] (.7*3,0) node {$\bullet$} ;
\draw[red] (.7*5,0) node {$\bullet$} ;

\draw[blue] (.7*2,0) node {$\times$} ;
\draw[blue] (.7*4,0) node {$\times$} ;
\draw[blue] (.7*6,0) node {$\times$} ;

\draw (.7*4,0) arc  (0:180:.7* .5) ;
\draw (.7*5,0) arc  (0:180:.7*1.5) ;

\draw (0+.7*8,0) node {\Huge $<$} ;

\draw (0+.7*9,0) --   (.7*7+.7*9,0) ;

\draw[red] (.7*1+.7*9,0) node {$\bullet$} ;
\draw[red] (.7*3+.7*9,0) node {$\bullet$} ;
\draw[red] (.7*5+.7*9,0) node {$\bullet$} ;

\draw[blue] (.7*2+.7*9,0) node {$\times$} ;
\draw[blue] (.7*4+.7*9,0) node {$\times$} ;
\draw[blue] (.7*6+.7*9,0) node {$\times$} ;

\draw (.7*3+.7*9,0) arc  (0:180:.7* .5) ;
\draw (.7*5+.7*9,0) arc  (0:180:.7* .5) ;

\draw (.7*8,-1.5) node {\Huge $\Downarrow$} ;

\draw (0+.7*4.5,-3) --   (.7*7+.7*4.5,-3) ;

\draw[red] (.7*1+.7*4.5,-3) node {$\bullet$} ;
\draw[red] (.7*3+.7*4.5,-3) node {$\bullet$} ;
\draw[red] (.7*5+.7*4.5,-3) node {$\bullet$} ;

\draw[blue] (.7*2+.7*4.5,-3) node {$\times$} ;
\draw[blue] (.7*4+.7*4.5,-3) node {$\times$} ;
\draw[blue] (.7*6+.7*4.5,-3) node {$\times$} ;

\draw (.7*4+.7*4.5,-3) arc  (0:180:.7* .5) ;
\end{tikzpicture}
\caption{Schematic representation of the result of Theorem~\ref{ThLMI}
  in the case $k_0=1$.}\label{pic:ThLMI}
\end{center}
\end{figure}

This result is represented in broad outline in Figure~\ref{pic:ThLMI}. 
For practical purposes,
these indicators allow to find pairs of neighbors that are
matched in an optimal transport plan. It also shows that the
usual $c$-cyclical monotonicity condition of
optimality~(see \cite{Smith1992328} and~\cite{springerlink:10.1007/BF02392620}) can be improved in the concave case:
only specific subsets have to be tested to check the optimality of a
transport plan.

\subsection{Algorithm}
\label{sec:algo}
%In this section, 
We now derive from the previous theorem %Theorem~\ref{ThLMI}
a simple algorithm to compute an optimal transport
plan in the case of chains. 
For the sake of simplicity, we only consider the balanced
case. The unbalanced case can be treated in the same way.

The local matching indicators defined in
Definition~\ref{def:LMI} can be used recursively to compute optimal
transport plans for chains. The elementary step of this approach consists in finding a negative
indicator satisfying the hypothesis of Theorem~\ref{ThLMI} in the list of supplies and
demands. Once this step is achieved, the inner points involved in  
this indicator are matched as prescribed in Theorem~\ref{ThLMI} and
removed from the list. %The research of negative indicators 
%then restarts.

 As for the research, it is performed by means of a loop that
updates iteratively the pair $(k_0,i_0)$, following
the lexicographic sorting (over admissible pairs), as long as positive
indicators are found. In this way, the hypothesis of 
Theorem~\ref{ThLMI} are satisfied when a negative indicator is found.
At the beginning of the algorithm or when a negative indicator
is found, the set of admissible pairs is updated and the current pair
is set to $(1,1)$.

We denote by $\sigma^\star$ the map for which this minimum is attained.
\medskip
\begin{algorithm}\label{alg:3}
\begin{itemize}
  \item Set $\mathcal{P}=\{p_1,\dots,p_N,q_1,\dots,q_N\}$, $\ell^p=\{1,\dots,N\}$, $\ell^q=\{1,\dots,N\}$, and $k=1$;
  \item while $\mathcal{P} \neq \emptyset $ and  $k < N$ 
  \begin{enumerate}
    \item compute $I^p_k(i)$ and $I^q_k(i')$ for $i=1,\dots,N-k$ and $i'=1,\dots,N-k-1$; \label{costlystep}
    \item define $$\mathcal{I}_k^p=\{i_0, 1\leq i_0\leq N-k ,I^p_k(i_0)<0 
     \},$$ $$\mathcal{I}_k^q=\{i_0, 1\leq i_0\leq N-k-1,I^q_k(i_0)<0 \} ;$$ 
    \item  if $\mathcal{I}_k^p=\emptyset$ and $\mathcal{I}_k^q= \emptyset$, then set $k=k+1$;
    \item  else do \label{removingstep}
    \begin{itemize}
      \item for all $i_0$ in $\mathcal{I}_k^p$ and for $i=i_0  +1,\dots,i_0+k$, do
      \begin{itemize}
        \item define $\sigma^\star(\ell^p_i)=\ell^q_{i-1}$,
        \item remove $\{p_{\ell^p_i},q_{\ell^q_{i-1}}\}$ from $\mathcal{P}$,
        \item remove $\ell^p_i$ and $\ell^q_i$ from $\ell^p$ and $\ell^q$ respectively;
      \end{itemize}
      \item for all $i_0'$ in $\mathcal{I}_k^q$ and for $i=i_0'+1,\dots,i_0'+k$, do
      \begin{itemize}
        \item define $\sigma^\star(\ell^p_i)=\ell^q_i$,
        \item remove $\{p_{\ell^p_i},q_{\ell^q_i}\}$ from $\mathcal{P}$,
        \item remove $\ell^p_i$ and $\ell^q_i$ from $\ell^p$ and $\ell^q$ respectively;
      \end{itemize}
      \item set $N=\frac 12 Card(\mathcal{P})$, and rename the points in
        $\mathcal{P}$ such that $$\mathcal{P}=\{p_1,\dots,p_N,q_1,\dots,q_N\},$$  
        $$p_1<q_1<\dots<p_{i}<q_{i}<p_{i+1}<q_{i+1}<\dots<p_{N}<q_{N}; $$
      \item set $k=1$;
    \end{itemize}
  \end{enumerate}
  \item if $k=N-1$, for $i=1,\dots,N$ set $\sigma^\star(\ell^p_i)=\ell^q_i$.
\end{itemize}
\end{algorithm}
\medskip
A first alternative algorithm consists in testing the sign of each
$I^p_k(i)$ and $I^q_k(i')$ as soon as they have been computed and
 in removing the corresponding pairs of points whenever a negative value is
found. A second alternative consists in carrying out the research of
pairs $(i_0,k_0)$ that satisfy the hypothesis of
Theorem~\ref{ThLMI} following the lexicographic order associated to the counter $(i_0+2k_0,k_0)$.
% What follows also holds for this variant. 

\subsection{Proof of Theorem~\ref{ThLMI}}
%\label{sec:theorem}

\subsubsection{Technical results}
%\label{sec:tech}
This section aims at introducing technical results that are
required to prove Theorem~\ref{ThLMI}. We keep the notations
introduced therein. We start with a basic result that plays a significant role in the
proof of Theorem~\ref{ThLMI}. As it was the case for the non-crossing
rule (Lemma~\ref{lemma:non-crossing-rule}), the  
concavity of the cost function is an essential assumption of this
lemma. 
\begin{lemma}\label{lemm:basic}%(Extension property)
We keep the previous notations. For $x,y \in \rset^{+}$, define 
$$
\varphi_{k,i}^p(x,y)=g(x+y+q_{i+k}-p_i)+\sum_{j=0}^{k-1}c(p_{i+j+1},q_{i+j})-g(x)-g(y)-\sum_{j=1}^{k-1}c(p_{i+j},q_{i+j}), 
$$
for $k,i\in \nset$, such that $1\leq k\leq N-1$ and $1\leq i\leq N-k$, and
$$
\varphi_{k,i}^q(x,y)=g(x+y+p_{i+k+1}-q_{i})+\sum_{j=1}^{k}c(p_{i+j},q_{i+j})-g(x)-g(y)-\sum_{j=1}^{k-1}c(p_{i+j+1}
,q_{i+j}),  
$$
for $k,i\in \nset$, such that $1\leq k\leq N-2$ and $1\leq i\leq N-k-1$ in the balanced case and $1\leq k\leq N-1$ and $1\leq i\leq N-k$ in the unbalanced case.
Both functions $
\varphi_{k,i}^p(x,y)$ and $
\varphi_{k,i}^q(x,y)$ are decreasing with respect to each of their two variables.
\end{lemma}

This lemma is a direct consequence of the concavity of the function
$g$. 
%\subsection{Tools specific to the unbalanced case}
%To deal with unbalanced chains, we need an additionnal result.

To deal with unbalanced chains, we need two additional lemmas, one of them requiring that $g$ is strictly monotone.
The first result we need is usually referred as ``The rule of three'' in the
literature~\cite{mccann1999est}.

\begin{lemma}[``rule of three'']
\label{lemma:ruleofthree}
Suppose that $g$ is strictly monotone. Given $p<q<p'<q'\in\rset$,
suppose that 
\begin{equation}\label{eq:r3}
c(p,q')+c(p',q) < c(p,q)+c(p',q').
\end{equation}
Then $ |p'-q|<\min(|p-q|,|p'-q'|)$.
\end{lemma}

\noindent{\em Proof:}
Since $g$ is increasing and since $|p-q'| \geq \max
(|p-q| ,|p'-q'| )$, Inequality~\eqref{eq:r3} implies that
$c(p',q) < \min (c(p,q) ,c(p',q') )$. The result follows the fact that
$g$ is strictly increasing.$\hfill \square\newline\newline$
We shall also make use of the following generalization. 
\begin{lemma}
\label{lemma:ruleofthreevariant}
Suppose that $g$ is strictly monotone.
  Under Hypothesis~(\ref{hyp1}) and (\ref{hyp3}) of Theorem~\ref{ThLMI}, the following inequalities are satisfied
  \begin{equation*}
%\label{generalized_ruleof3}
    |q_i-p_{i+1}|<\min(|p_{i_0}-q_i|,|p_{i+1}-q_{i_0+k_0}|), \;\;\forall i \in \{i_{0},\dots,i_0+k_0-1\}.
  \end{equation*}
If Hypothesis~\ref{hyp3} is replaced by $I^q_{k_0}(i_0')< 0$ and
Hypothesis~(\ref{hyp2}) holds, one finds
  \begin{equation*}
%\label{generalized_ruleof3_2}
    |p_i-q_i|<\min(|q_{i_0'}-p_i|,|q_i-p_{i_0'+k_0+1}|), \;\;\forall i \in \{i_{0}'+1,\dots,i_0'+k_0\}.
  \end{equation*}

\end{lemma}

%\begin{proof}
\noindent{\em Proof:} Let $i\in \{i_{0},\dots,i_0+k_0-1\}$. Hypothesis~\eqref{hyp3} of Theorem~\ref{ThLMI} implies that
$$
    c(p_{i+1},q_i) + c(p_{i_0},q_{i_0+k_0}) < \sum_{j=i_0}^{i_0+k_0}
    c(p_j,q_j) - \sum_{j=i_0}^{i_0+k_0-1} c(p_{j+1},q_j) +
    c(p_{i+1},q_i). 
$$
Now, because of Hypothesis~\eqref{hyp1}, we have
$I^p_{i-i_0}(i_0)\geq 0$ and $I^p_{i_0+k_0-i-1}(i+1)\geq 0$, which
means that 
$$\sum_{j=i_0}^{i}c(p_j,q_j) \le c(p_{i_0},q_i) + \sum_{j=i_0}^{i-1} c(p_{j+1},q_j)$$ and
$$\sum_{j=i+1}^{i_0+k_0}
c(p_j,q_j) \le c(p_{i+1},q_{i_0+k_0}) + \sum_{j=i+1}^{i_0+k_0-1}
c(p_{j+1},q_j).$$
Thus,  
$$
  c(p_{i+1},q_i) + c(p_{i_0},q_{i_0+k_0}) < c(p_{i_0},q_i) +  c(p_{i+1},q_{i_0+k_0}). 
$$
We conclude with the rule of three.
The result in the case $I^q_{k_0}(i_0')< 0$ can be deduced by
symmetry. $\hfill \square\newline\newline$

%\end{proof}

Note that in the two previous proofs, the only necessary hypothesis is
that the cost is a strictly increasing
function of the distance. In particular the result also
holds in the case where the cost
function is increasing and convex.

\begin{lemma}[``partial sums'']
\label{lemma:partialsums}
 Under Hypothesis~(\ref{hyp1})  and~(\ref{hyp3}) of Theorem~\ref{ThLMI}, for any $i$ in
 $\{i_0+1,\dots,i_0+k_0\}$ and $i'$ in $\{i_0,\dots,i_0+k_0-1\}$, the
 following inequalities are satisfied:
\begin{equation}
  \label{eq:startsum}
\sum_{j=i_0}^{i-1} c(p_{j},q_j)  > \sum_{j=i_0}^{i-1} c(p_{j+1},q_j),
\end{equation}
and
\begin{equation}
  \label{eq:endsum}
\sum_{j=i'+1}^{i_0+k_0} c(p_{j},q_j)  > \sum_{j=i'}^{i_0+k_0-1} c(p_{j+1},q_j).
\end{equation}
If Hypothesis~\ref{hyp3} is replaced by $I^q_{k_0}(i_0')< 0$ and
Hypothesis~(\ref{hyp2}) holds, one finds
\begin{equation}
  \label{eq:startsum2}
\sum_{j=i_0}^{i-1} c(p_{j+1},q_j)  > \sum_{j=i_0+1}^{i} c(p_j,q_j),
\end{equation}
and
\begin{equation}
  \label{eq:endsum2}
\sum_{j=i'+1}^{i_0+k_0} c(p_{j+1},q_j)  > \sum_{j=i'+1}^{i_0+k_0} c(p_{j},q_j).
\end{equation}
\end{lemma}
%\begin{proofof}

\noindent{\em Proof:} In order to prove
inequality~\eqref{eq:startsum}, remark that since $I^p_{i_0}(k_0)<0$
\begin{eqnarray*}
  \sum_{j=i_0}^{i-1} c(p_j,q_j) &=&  \sum_{j=i_0}^{i_0+k_0} c(p_j,q_j)
  -  \sum_{j=i}^{i_0+k_0} c(p_j,q_j)\\ 
&>& c(p_{i_0},q_{i_0+k_0})+\sum_{j=i_0}^{i_0+k_0-1} c(p_{j+1},q_j) -
\sum_{j=i}^{i_0+k_0} c(p_j,q_j),
\end{eqnarray*}
for $i$ such that $i_0+1\leq i\leq i_0+k_0$.
Moreover, since $I^p_{i_0+k_0-i}(i) \geq 0$ one has
$$
\sum_{j=i_0}^{i-1} c(p_j,q_j) >
c(p_{i_0},q_{i_0+k_0})+\sum_{j=i_0}^{i_0+k_0-1} c(p_{j+1},q_j) -
c(p_{i},q_{i_0+k_0}) -  \sum_{j=i}^{i_0+k_0-1} c(p_{j+1},q_j).
$$
Since $g$ is increasing, this leads to the inequality~\eqref{eq:startsum}. The proof of
Equations~(\ref{eq:endsum}--\ref{eq:endsum2}) follows the same path. $\hfill \square\newline\newline$

We are now in the position to prove our main result. In a first part
we focus on the balanced case, and then go to the unbalanced case,
which requires more efforts.

\subsubsection{The balanced case}\label{sec:proofbal}Consider the balanced case,
i.e., the situation corresponding to \eqref{balanced}.
We focus on the case where $I^p_{k_0}(i_0)< 0$. The case
$I^q_{k_0}(i_0')<0$ can be treated the same way.

The proof consists in proving that Hypothesis~(\ref{hyp1}--\ref{hyp3}) of Theorem~\ref{ThLMI} imply that
neither demand nor supply points located between $p_{i_0}$ and
$q_{i_0+k_0}$ %%%AS: $q_{i_0+k_0}$ ? JS : done.
can be matched with points located outside this interval,
i.e. that the set $\mathcal{S}^{k_0}_{i_0}=\{ p_j, i_0+1\leq j \leq i_0+k_0 \}\cup \{q_j,
i_0\leq j \leq i_0+k_0-1\}$ is invariant under an optimal transport plan. In
this case, the result follows from Hypothesis~(\ref{hyp1}--\ref{hyp2}).

Suppose that $\mathcal{S}^{k_0}_{i_0}$ is not preserved by an optimal
transport plan $\sigma^\star$. According to the non-crossing rule, three cases can occur:
\begin{enumerate}[a)]
\item\label{c1} There exists $i_1\in \nset$, such that $1\leq i_1\leq i_0$ and
  $ i_0 \leq \sigma^\star(i_1)\leq i_0+k_0 -1$ and there exists $i_1'\in
  \nset$, such that $\sigma^\star(i_1)+1 \leq i_1'\leq i_0+k_0$ and
  $ i_0+k_0 \leq \sigma^\star(i_1')\leq N$.
\item\label{c2} There exists $i_2\in \nset$, with $i_0+1\leq i_2\leq i_0+k_0$ such that
  $ 1 \leq \sigma^\star(i_2)\leq i_0-1$.
\item\label{c3} There exists $i_2\in
  \nset$, with $i_0+k_0 < i_2 \leq  N$ such that
  $ i_0 \leq \sigma^\star(i_2)< i_0+k_0$.
\end{enumerate}
We first prove that Case~\ref{c1}) cannot occur.

In Case~\ref{c1}), one can assume without loss of generality that $\sigma^\star (i_1)$ is the largest
index such that $1\leq i_1\leq i_0$, $ i_0 \leq \sigma^\star(i_1)\leq
i_0+k_0 -1$ and that $i_1'$ is the smallest index such that
$\sigma^\star(i_1)+1 \leq i_1'\leq i_0+k_0$, $ i_0+k_0 \leq \sigma^\star(i_1')\leq N$ 
. Assume also that we are not in Cases~\ref{c2})  or~\ref{c3}). With
such assumptions and because of the non-crossing rule, the (possibly
empty) subset $\{ p_i, \sigma^\star (i_1)+1 \leq i \leq i_1'-1 \}\cup
\{q_i, 
\sigma^\star (i_1)+1\leq i \leq i_1'-1\}$ is stable by
$\sigma^\star$. Because of Hypothesis~(\ref{hyp1}--\ref{hyp2}), no 
nesting ({\em i.e.} no pair of nested matchings) can occur in this subset, and $\sigma^\star(i)=i$ for
$i=\sigma^\star(i_1)+1,\dots,i_1'-1$.

On the other hand, since $\sigma^\star$ is optimal, one has:
$$
c(p_{i_1},q_{\sigma^\star(i_1)})+c(p_{i_1'},q_{\sigma^\star(i_1')})+\sum_{j=\sigma^\star(i_1)+1}^{i_1'-1}c(p_j,q_j)
\leq
c(p_{i_1},q_{\sigma^\star(i_1')})+\sum_{j=\sigma^\star(i_1)}^{i_1'-1}c(p_{j+1},q_j).
$$
Thanks to Lemma~\ref{lemm:basic},
one deduces from this last inequality that:
$$
c(p_{i_0},q_{\sigma^\star(i_1)})+c(p_{i_1'},q_{i_0+k_0})+\sum_{j=\sigma^\star(i_1)+1}^{i_1'-1}c(p_j,q_j)\\
\leq
c(p_{i_0},q_{i_0+k_0})+\sum_{j=\sigma^\star(i_1)}^{i_1'-1}c(p_{j+1},q_j),
$$
and then:
\begin{eqnarray}
c(p_{i_0},q_{\sigma^\star(i_1)})+\sum_{j=i_0}^{\sigma^\star(i_1)-1}c(p_{j+1},q_j)
+c(p_{i_1'},q_{i_0+k_0})+\sum_{j=i_1'}^{i_0+k_0-1}c(p_{j+1},q_j)
\nonumber\\
+\sum_{j=\sigma^\star(i_1)+1}^{i_1'-1}c(p_j,q_j)\leq
c(p_{i_0},q_{i_0+k_0})+\sum_{j=i_0}^{i_0+k_0-1}c(p_{j+1},q_{j}).\ \ \label{ineq:1}
\end{eqnarray}
According to Hypothesis~(\ref{hyp1}),
$I^p_{\sigma^\star(i_1)-i_0}(i_0)\geq 0$ and
$I^p_{i_0+k_0-i_1'}(i_1')\geq 0$, so that:

$$ \sum_{j=i_0}^{\sigma^\star(i_1)}c(p_j,q_j)\leq
c(p_{i_0},q_{\sigma^\star(i_1)})+\sum_{j=i_0}^{\sigma^\star(i_1)-1}c(p_{j+1},q_{j})$$
$$\sum_{j=i_1'}^{i_0+k_0} c(p_j,q_j)\leq
c(p_{i_1'},q_{i_0+k_0})+\sum_{j=i_1'}^{i_0+k_0-1}c(p_{j+1},q_j).
$$
Combining these last inequalities with~(\ref{ineq:1}) one finds that:
$$\sum_{j=i_0}^{i_0+k_0}c(p_j,q_j) \leq c(p_{i_0},q_{i_0+k_0})+\sum_{j=i_0}^{i_0+k_0-1}c(p_{j+1},q_j),
$$
which contradicts Hypothesis~(\ref{hyp3}). 

Let us now prove that Cases~\ref{c2}) and~\ref{c3}) contradict the
assumptions. 
As Cases~\ref{c2}) and~\ref{c3}) can be treated in the same
way, we only consider Case~\ref{c2}). Without loss of generality, one can
assume that $i_2$ is the smallest index such that $i_0+1\leq i_2 \leq
i_0+k_0$ and $\sigma^\star(i_2)\leq i_0-1$. Because of the
non-crossing rule and the fact there are
necessarily as many demands as supplies between $q_{i_0}$ and
$p_{i_2}$, there exists one and only one index $i_2'$ such that
$i_0\leq \sigma^\star(i_2')\leq i_2-1$ and $1\leq i_2'\leq
i_0$. Consequently, the non-crossing rule implies that the (possibly
empty) subsets $\{ p_i, i_0 +1 \leq i \leq \sigma^\star(i_2')  \}\cup
\{q_i, 
i_0   \leq i \leq \sigma^\star(i_2')-1\}$ and $\{ p_i,
\sigma^\star(i_2')+1 \leq i \leq i_2-1 \}\cup \{q_i,
\sigma^\star(i_2')+1\leq i \leq i_2 -1\}$ are stable by an optimal
transport plan. Because of Hypothesis~(\ref{hyp1}--\ref{hyp2}), no
nesting can occur in these subsets, and $\sigma^\star(i)=i-1$ for
$i=i_0+1,\dots,\sigma^\star(i_2')$ and $\sigma^\star(i)=i$ for
$i=\sigma^\star(i_2')+1,\dots,i_2-1$. 

On the other hand, since $\sigma^\star$ is optimal, one has
\begin{eqnarray*}
c(p_{i_2},q_{\sigma^\star(i_2)})+c(p_{i_2'},q_{\sigma^\star(i_2')})+\sum_{j=i_0+1}^{\sigma^\star(i_2')}c(p_{j},q_{j-1})
+\sum_{j=\sigma^\star(i_2')+1}^{i_2-1}c(p_{j},q_{j})
\\\leq c(p_{i_2'},q_{\sigma^\star(i_2)}) + \sum_{j=i_0+1}^{i_2}c(p_{j},q_{j-1}).
\end{eqnarray*}
Thanks to Lemma~\ref{lemm:basic},
one deduces from this last inequality that:
\begin{eqnarray}
c(p_{i_2},q_{\sigma^\star(i_2)})+c(p_{i_0},q_{\sigma^\star(i_2')})+\sum_{j=i_0+1}^{\sigma^\star(i_2')}c(p_{j},q_{j-1})
+\sum_{j=\sigma^\star(i_2')+1}^{i_2-1}c(p_{j},q_{j})\nonumber
\\\leq c(p_{i_0},q_{\sigma^\star(i_2)}) + \sum_{j=i_0+1}^{i_2}c(p_{j},q_{j-1}).\label{ineq3}
\end{eqnarray}
Because the cost is supposed to be increasing with respect to the
distance, one finds that
$c(p_{i_0},q_{\sigma^\star(i_2)})\leq c(p_{i_2},q_{\sigma(i_2)})$, so that~(\ref{ineq3}) implies:
%\begin{eqnarray*}
$$c(p_{i_0},q_{\sigma^\star(i_2')})+\sum_{j=i_0+1}^{\sigma^\star(i_2')}c(p_{j},q_{j-1})
+\sum_{j=\sigma^\star(i_2')+1}^{i_2-1}c(p_{j},q_{j})
\leq \sum_{j=i_0+1}^{i_2}c(p_{j},q_{j-1}),$$
%\end{eqnarray*}
and then:
\begin{eqnarray}
c(p_{i_0},q_{\sigma^\star(i_2')})+\sum_{j=i_0+1}^{\sigma^\star(i_2')}c(p_{j},q_{j-1})
+\sum_{j=\sigma^\star(i_2')+1}^{i_2-1}c(p_{j},q_{j})+\sum_{j=i_2+1}^{i_0+k_0}c(p_{j},q_{j-1})\nonumber\\
\leq \sum_{j=i_0+1}^{i_0+k_0}c(p_{j},q_{j-1}).\label{ineq4}
\end{eqnarray}
According to Hypothesis~(\ref{hyp1}), $I^p_{\sigma^\star(i_2')-i_0}(i_0)\geq 0$, so that:

$$\sum_{j=i_0}^{\sigma^\star(i_2')}c(p_j,q_j)\leq
c(p_{i_0},q_{\sigma^\star(i_2')})+\sum_{j=i_0}^{\sigma^\star(i_2')-1}c(p_{j+1},q_j).$$
Combining these last inequalities with~\eqref{ineq4} one finds that:
$$\sum_{j=i_0}^{i_0+k_0}c(p_j,q_j) \leq c(p_{i_0},q_{i_0+k_0})+\sum_{j=i_0}^{i_0+k_0-1}c(p_{j+1},q_j),
$$
which contradicts Hypothesis~(\ref{hyp3}).

We have then shown that 
neither demand nor supply points located between $p_{i_0}$ and
$q_{i_0+k_0+1}$ can be matched with points located outside this interval. The
set $\mathcal{S}^{k_0}_{i_0}$ is then stable by an optimal transport
plan. According to Hypothesis~(\ref{hyp1}--\ref{hyp2}), no nesting can
occur in $\mathcal{S}^{k_0}_{i_0}$. The result follows.  $\hfill \square\newline\newline$

\subsubsection{The unbalanced case}
\label{sec:unbalanced}
%In this section, 
We then  %assume that $g$ is strictly monotone and
show that Theorem~\ref{ThLMI} still holds in the
unbalanced case. We start with the case $I^p_{k_0}(i_0) < 0$.\\
%In this section, we show that Theorem~\ref{ThLMI} still holds in the unbalanced case.
Observe first that none of the points $p_{j}$,
$i_0+1\leq j \leq i_0+k_0$ can remain unmatched in an optimal
transport plan.
Indeed, assume on the contrary that there exists $\ell$ in $\{i_0+1,\dots, i_0+k_0\}$ such
that $p_\ell$ is unmatched in an optimal transport plan
$\sigma^\star$.   
Note first that no nesting can occur in $\mathcal{S}^{k_0}_{i_0}$, so that
the points in this set can only be matched either with their neighbors or with points outside this set.
According to Lemma~\ref{lemma:partialsums} 
$$
\sum_{j = i_0}^{\ell - 1} c(p_j, q_j)
> \sum_{j = i_0}^{\ell - 1} c(p_{j + 1}, q_j).
$$
Therefore  
%and since no nesting can occur in $\mathcal{S}^{k_0}_{i_0}$, 
we cannot have
$\sigma^\star(i)=i$ for $i=i_0,\dots,\ell-1$: otherwise it would be 
possible to rematch all the points $q_i$ in this interval to their
right neighbors and reduce the cost. Hence, as the point~$p_\ell$ is
unmatched and, because of Lemma~\ref{lemma:unma-isol}, exposed, there exists $m$ in
$\{i_0,\dots, \ell-1\}$ such that $(\sigma^\star)^{-1}(m) <i_0$. 
Choose $m$ to be the greatest value of the index satisfying this property.
 Since no nesting can occur in $\mathcal{S}^{k_0}_{i_0}$, we have $\sigma^\star(i) = i$ for all~$i$ in the
(possibly empty) interval $m + 1\le i\le \ell - 1$. 
Now, since $g$ is an increasing function,
$$
  c(p_{(\sigma^\star)^{-1}(m)},q_m) + \sum_{j=m+1}^{\ell-1} c(p_j,q_j) 
> c(p_{i_0},q_m) +\sum_{j=i_0}^{\ell-1} c(p_j,q_j) - \sum_{j=i_0}^{m} c(p_j,q_j)
$$
Using again Equation~\eqref{eq:startsum} of Lemma~\ref{lemma:partialsums}, one deduces from this last
inequality that
$$
  c(p_{(\sigma^\star)^{-1}(m)},q_m) + \sum_{j=m+1}^{\ell-1} c(p_j,q_j) > c(p_{i_0},q_m)
+\sum_{j=i_0}^{\ell-1} c(p_{j+1},q_j) - \sum_{j=i_0}^{m} c(p_j,q_j).$$
It follows from this and from $I^p_{m-i_0}(i_0) \geq 0$ that
\begin{eqnarray*}
  c(p_{(\sigma^\star)^{-1}(m)},q_m) + \sum_{j=m+1}^{\ell-1} c(p_j,q_j) &>&
  c(p_{i_0},q_m) +\sum_{j=i_0}^{m-1}
  c(p_{j+1},q_j)\\&&+\sum_{j=m}^{\ell-1} 
c(p_{j+1},q_j) - \sum_{j=i_0}^{m} c(p_j,q_j)\\ 
&\geq &\sum_{j=m}^{\ell-1} c(p_{j+1},q_j).
\end{eqnarray*}
In other words, it is cheaper to match each $q_i$, $m\le i\le \ell -
1$, with its right neighbor $p_{i+1}$ and to exclude
$p_{(\sigma^\star)^{-1}(m)}$ than to match each 
$q_i$ with its neighbor $p_{i}$ and to exclude $p_{\ell}$. In all cases,
the point $p_\ell$ cannot remain unmatched. 

If the point $p_{i_0}$ is matched in the transport plan $\sigma^\star$,
%and if $\sigma^\star$ is optimal, 
then we can conclude by the already
proved first part of Theorem~\ref{ThLMI} that $\sigma^\star(i)=i-1$
for $i=i_0+1,\dots,i_0+k_0$ %if $\sigma^\star$ is optimal 
(according to Lemma~\ref{lemma:unma-isol}  
unmatched points are exposed, the existence of an unmatched $p_i$
outside of $[p_{i_0},q_{i_0+k_0}]$ has no consequence on this result).

Now, assume that $p_{i_0}$ remains unmatched and that there exists $m'$
in $\{i_0,\dots,i_0+k_0-1\}$ such that $(\sigma^\star)^{-1}(m')\neq m'+1$. Since
$p_{i_0}$ is exposed,  and
since all points of ${\cal S}_{i_0}$ are matched and no nesting
can occur in ${\cal S}_{i_0}$, there exists $m'$  in $\{ i_0,\dots, i_0
+k_0 -1\}$ such that $(\sigma^\star)^{-1} (m') >
i_0 + k_0$.
% since matched points are either neighbors
% or separated by more than $2k_0$ points, $(\sigma^\star)^{-1}(m)> i_0 +k_0$.
One can assume without loss of generality that $m'$ is the largest
index in $\{i_0,\dots,i_0+k_0-1\}$ satisfying $(\sigma^\star)^{-1}(m')> i_0
+k_0$. 

Actually, $m'<i_0+k_0-1$. Indeed, suppose that $m'=i_0+k_0-1$:
on the one hand, because of Hypothesis~(\ref{hyp1}--\ref{hyp3}), the
rule of three (variant, Lemma~\ref{lemma:ruleofthreevariant}) implies that
$|p_{i_0+k_0}-q_{i_0+k_0-1}|<|p_{i_0+k_0}-q_{i_0+k_0}|$. But on the
other hand, since the matchings $(p_{(\sigma^\star)^{-1}(m')},q_{m'})$ and
$(p_{i_0+k_0},q_{\sigma^\star(i_0+k_0)})$ belong to an optimal
transport plan, the rule of three (standard version,
Lemma~\ref{lemma:ruleofthree}) implies
$|p_{i_0+k_0}-q_{i_0+k_0-1}|>|p_{i_0+k_0}-q_{\sigma^\star(i_0+k_0)}|$. Because
of the non-crossing rule, $\sigma^\star(i_0+k_0)\geq i_0+k_0 $, hence 
$|p_{i_0+k_0}-q_{i_0+k_0-1}|>|p_{i_0+k_0}-q_{i_0+k_0}|$. This provides a contradiction.
%$(\sigma^\star)^{-1}(i_0+k_0-1) \leq i_0 +k_0$, otherwise there would exist
%$q_j$, with $j\geq i_0+k_0$ such that
%$|p_{i_0+k_0}-q_{j}|<|p_{i_0+k_0}-q_{i_0+k_0-1}|$ (again, because of
%the rule of three), which contradicts the previous inequality. It
%follows that $m<i_0+k_0-1$. 

Two cases can now occur: either  $\sigma^\star(i) =
i$ for all $i$ in $\{m'+1,\dots, i_0+k_0\}$, or there exists a unique supply $p_k$ in $\{m'+1,\dots,i_0+k_0\}$ such that
$\sigma^\star(k)>i_0+k_0$. This cannot happen for two different supplies in
$\{m'+1,\dots,i_0+k_0\}$, otherwise there would be another demand $q_\ell$
between these supplies such that $(\sigma^\star)^{-1}(\ell)>i_0+k_0$. 

In the first case, thanks to equation~\eqref{eq:endsum}
\begin{eqnarray*}
c(q_{m'} , p_{(\sigma^\star)^{-1}(m')})+ \sum_{j=m'+1}^{i_0+k_0}
c(p_j,q_j)&>& c(q_{m'} , p_{(\sigma^\star)^{-1}(m')})+
\sum_{j=m'}^{i_0+k_0-1} c(p_{j+1},q_j) \\ 
&>& c(q_{i_0+k_0},p_{(\sigma^\star)^{-1}(m')})+ \sum_{j=m'}^{i_0+k_0-1} c(p_{j+1},q_j),
\end{eqnarray*}
which contradicts the optimality of $\sigma^\star$.

In the second case, since $I^p_{k_0}(i_0)<0$,
\begin{eqnarray*}
c(q_{m'},p_{(\sigma^\star)^{-1}(m')}) + \sum_{j=m'+1}^{k-1} c(p_j,q_j)
+c(p_k,q_{\sigma^\star(k)}) > c(q_{m'},p_{(\sigma^\star)^{-1}(m')}) + c(p_k,q_{\sigma^\star(k)})\\ +
c(p_{i_0},q_{i_0+k_0}) +\sum_{j=i_0}^{i_0+k_0-1} c(p_{j+1},q_j) -
\sum_{j=i_0}^{m'} c(p_{j},q_j)-\sum_{j=k}^{i_0+k_0} c(p_{j},q_j). 
\end{eqnarray*}
Now, since $I^p_{m'-i_0}(i_0)\geq 0$ and $I^p_{i_0+k_0-k}(k)\geq 0$, this inequality yields
\begin{eqnarray*}
c(q_{m'},p_{(\sigma^\star)^{-1}(m')}) + \sum_{j=m'+1}^{k-1} c(p_j,q_j) +c(p_k,q_{\sigma^\star(k)})
% > c(q_m,p_{(\sigma^\star)^{-1}(m)}) + c(p_k,q_{\sigma^\star(k)}) +
% c(p_{i_0},q_{i_0+k_0})+ &\sum_{i=i_0}^{i_0+k_0-1} c(p_{i+1},q_i)& -
% c(p_{i_0},q_m) - c(p_k,q_{i_0+k_0}) - \sum_{i=i_0}^{m-1}
% c(p_{i+1},q_i)-\sum_{i=k}^{i_0+k_0-1} c(p_{i+1},q_i)\\ 
>  c(q_{m'},p_{(\sigma^\star)^{-1}(m')}) \\+ c(p_k,q_{\sigma^\star(k)}) -
  c(p_k,q_{i_0+k_0}) + c(p_{i_0},q_{i_0+k_0})
  - c(p_{i_0},q_{m'}) + \sum_{j=m'}^{k-1} c(p_{j+1},q_j).
\end{eqnarray*}
The two differences that appear in the right-hand side are positive so that
\begin{eqnarray*}
c(q_{m'},p_{(\sigma^\star)^{-1}(m')}) + \sum_{j=m'+1}^{k-1} c(p_j,q_j) +c(p_k,q_{\sigma^\star(k)})
\geq c(q_{\sigma^\star(k)},p_{(\sigma^\star)^{-1}(m')})\\ + \sum_{j=m'}^{k-1}
c(p_{j+1},q_j),
\end{eqnarray*}
which also contradicts the optimality of $\sigma^\star$.

By symmetry, the theorem remains valid in the case where  $I^q_{k_0}(i'_0)<0$
instead of ${I^p_{k_0}(i_0)<0}$.$\hfill \square\newline\newline$
%\end{proofof}

%--------------------------------------------------------------------
%               SECTION IV
%--------------------------------------------------------------------
\section{General unitary case}\label{sec:unitary}

% Note that a fully used supply~$p_{i_0}$ may too be exposed in a
% transport plan~$\gamma$ if all intervals of the form
% $(\min(p_{i_0},q_{j_0}), \max(p_{i_0},q_{j_0}))$ such that
% $\gamma_{i_0 j_0} > 0$ do not intersect with intervals corresponding
% to other arcs with positive~$\gamma_{ij}$, $i\neq i_0$.

We now focus on the general \emph{unitary} case, i.e. situations where $s_i = d_j =
1$ for all $i$,~$j$ and therefore $S = M$ and $D = N \le M$. 
As a main result, we shall explain how this case can be recast in
independent problems involving chains.
Recall that as it was the case for chains, in such a framework, optimal
transport plans are described by permutations.
%problems
%are actually assignment problems, where masses cannot be cut.

%\subsection{local balance of supplies and demands}
A consequence of the non-crossing rule~\ref{lemma:non-crossing-rule}  
 is usually called the {\it local balance
  of supplies and demands}: in the unitary case,
there are as many supplies as demands between any two
matched points $p_{i_0}$ and $q_{j_0}$.
% whereas for general real
%supplies and demands, the total supply and the total demand within
%the interval corresponding to $p_{i_0}$ and~$q_{j_0}$ with
%$\gamma_{i_0j_0} > 0$ may differ but balance can always be achieved 
%by including suitable shares of~$s_{i_0}$ and~$d_{j_0}$ (for details
%see subsection~\ref{sec:chains-real-valued}).
% Therefore, unmatched supplies are isolated
% on the line in the sense that whenever $\sigma^\star(i_0)=j_0$, no unmatched supply
% $p_{i}$ can be located between $p_{i_0}$ and $q_{j_0}$.
% defined in the
%following two subsections. 
%\subsection{Chains in the unitary case}
%\label{sec:chains-unitary-case}
% In this section, we present a way to subdivide the initial set of
% points $P\cup Q$ into a family of particular subsets called {\it
%   Chains} that are preserved by optimal transport plans. As a result,
% we obtain a precomputation step that enables 
% to break down an optimal transport problem into a set of smaller and
% independent problems.
We derive from this property a definition of chains in the case of unit masses.
Given a supply point $p_i$, define its {\it left
neighbor} $q'_i$ as the nearest demand point on the left of $p_i$ such that
the numbers of supplies and demands in the interval $(q'_i, p_i)$ are equal;
define the {\it right neighbor} $q''_i$ of $p_i$ in a similar
way.
Furthermore define left and right neighbors of a demand point~$q_j$ to
be the supply points that have~$q_j$ as  their right and
left neighbor, respectively.
Iterating this procedure gives raise to {\it a
  chain}.
\newcounter{tempi}
\setcounter{tempi}{\value{theorem}}
\begin{definition}[unitary case]
\label{def:chain-int}
A \emph{chain} in $P\cup Q$ is a maximal alternating sequence of
supplies and demands of one of the forms
\begin{enumerate}
\item $(p_{i_1}, q_{j_1}, \dots, p_{i_k}, q_{j_k})$,\label{cas1}
\item $(q_{j_1}, p_{i_1}, \dots, q_{j_k}, p_{i_k})$, \label{cas2}
\item $(p_{i_1}, q_{j_1}, \dots, q_{j_{k-1}}, p_{i_k})$,\label{cas3}
\end{enumerate}
with $k\geq 1$ and such that each pair of consecutive points in
the sequence is made of a point and its right neighbor.
\end{definition}

\begin{figure}
  \begin{center}
    \begin{tikzpicture}
      \draw (0,0) -- (9,0) ;

      \draw[red] (1,0) node {$\bullet$} (4,0) node {$\bullet$} %
      (5,0) node {$\bullet$} (8,0) node {$\bullet$} ;

      \draw[blue] (2,0) node {$\times$} (3,0) node {$\times$} %
      (6,0) node {$\times$} (7,0) node {$\times$} ;

      \draw[dashed] (2,0) arc (0:180: .5) %
      (5,0) arc (0:180:1.5) %
      (6,0) arc (0:180: .5) ;

      \draw[dashed] (4,0) arc (0:-180: .5) %
      (7,0) arc (0:-180:1.5) %
      (8,0) arc (0:-180: .5) ;
    \end{tikzpicture} \\[1ex]
    \begin{tikzpicture}
      \draw (0, 0) -- (9, 0) ;

      \draw (1, 1) -- (2, 1) (3, -1) -- (4, -1) %
      (5, 1) -- (6, 1) (7, -1) -- (8, -1) ;

      \draw[dashed] (1, .5) -- (6, .5) %
      (3, -.5) -- (8, -.5) ;

      \draw[red] (1, 0) node {$\bullet$} -- (1, 1) %
      (4, -1) node {$\bullet$} -- (4, 0) %
      (5, 0) node {$\bullet$} -- (5, 1) %
      (8, -1) node {$\bullet$} -- (8, 0) ;

      \draw[blue] (2, 1) node {$\times$} -- (2, 0) %
      (3, 0) node {$\times$} -- (3, -1) %
      (6, 1) node {$\times$} -- (6, 0) %
      (7, 0) node {$\times$} -- (7, -1) ;
    \end{tikzpicture}
 \caption{\label{exchain}  Example of a problem containing two chains. Top: chains   represented as collections of dashed arcs. Bottom: chains represented as dashed lines connecting elements of mass that are   left and right neighbors (cf Figure~\ref{fig:chainexnonint}).}
  \end{center}
\end{figure}
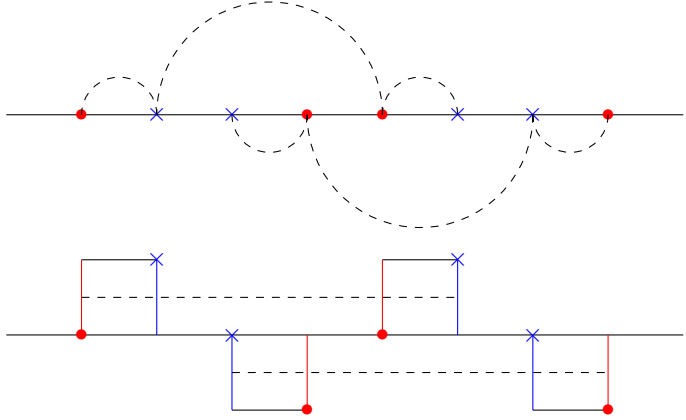

Examples of chains are shown in Figure~\ref{exchain}.
Observe that because of Case~(\ref{cas3}), some chains can be composed
of only one (unmatched) supply, and no
demand.

Because of the local balance property, matching in an optimal
plan can occur only between points that belong to the same chain.
Indeed, an extension of the
proof of Lemma 3 of \cite{aggarwal1992efficient} %%%AS: add a detailed
                                %%%proof? JS : I have added a
                                %%%sentence. I have the feeling that
                                %%%it is just a consequence of the stability. 
shows that the family of chains forms a
partition of $P \cup Q$.  As a consequence, each chain is preserved
by an optimal transport plan. For example, if a chain is composed of a single supply, it cannot
be matched in any optimal transport plan and can thus be dismissed
from the problem. In summary, the general
unitary problem can be decomposed into independent problems 
that only deal with chains, and then apply the results of Section~\ref{sec:chains}.

Note finally that the construction of the
set of chains only depends on relative positions of supplies and demands
and does not involve any evaluation of the cost function.
The exact construction is not described here because it
is subsumed by the algorithm presented in
subsection~\ref{sec:data-struct-algor}.
% It can be achieved in $O(N+M)$ operations:
% a possible  algorithm  consists in considering sequentially (say from
% the left to the right) the
% points of the set $P\cup Q$, and, at the same time,  building all the chains
% iteratively either by adding the current point to an existing chains,
% or by initiating a new chain in such a way that the local balance
%   of supplies and demands is fulfilled.  
% Note also that this partition makes it possible to parallelize
% transport optimization algorithms.

\bigskip

%--------------------------------------------------------------------
%               SECTION IV
%--------------------------------------------------------------------
\section{Non-unitary case}\label{sec:nonunitary}
\subsection{Chains in real-valued histograms}
\label{sec:chains-real-valued}

In this case the notions of right and left neighbors should be defined
for \emph{infinitesimal elements} of supply and demand.  The
corresponding definition may be given in purely intrinsic terms, but
the following graphical representation makes it more evident.

Consider the signed measure $\sum_i s_i \delta_{p_i} - \sum_j d_j
\delta_{q_j}$ on the real line, where $\delta_x$ is a unit Dirac mass
at~$x$.  Plot its cumulative distribution function~$F$, whose graph
has an upward jump at each~$p_i$ and a downward jumps at each~$q_j$,
and augment it with vertical segments to make the graph into a
continuous curve (Figure~\ref{fig:chainexnonint}; cf also Figure~\ref{exchain}).  Thus, e.g., the
segment corresponding to a supply point~$p_i$ connects the points of
the graph with coordinates $(p_i, F(p_i))$ and~$(p_i, F(p_i + 0) =
F(p_i) + s_i)$ (assuming left continuity of~$F$).  Here and below in
figures similar to Figure~\ref{fig:chainexnonint} vertical segments
corresponding to supply points are plotted in red and those
corresponding to demand points in blue (color online).

\begin{figure}
  \centering
  \begin{tikzpicture}
    \draw (-.5, 0) -- (8.5, 0) ;

    \draw (0, 4.75) -- (1.5, 4.75) %
    (1.5, 1.5) -- (2.25, 1.5) %
    (2.25, 5.25) -- (3.25, 5.25) %
    (3.25, .5) -- (4.5, .5) %
    (4.5, -1.25) -- (5.75, -1.25) %
    (5.75, 2.5) -- (6.5, 2.5) %
    (6.5, -.5) -- (8, -.5) %
    (8, 1) -- (8.5, 1);

    \draw[dotted] %
    (-.5, 5.25) node[left] {$y_0 = F(q_2)$} -- (8.5, 5.25) %
    (-.5, 4.75) node[left] {$y_1 = F(q_1)$} -- (8.5, 4.75) %
    (-.5, 2.5) node[left] {$y_2 = F(q_4)$} -- (8.5, 2.5) %
    (-.5, 1.5) node[left] {$y_3 = F(p_2)$} -- (8.5, 1.5) %
    (-.5, 1) node[left] {$y_4 = F(p_4 + 0)$} -- (8.5, 1) %
    (-.5, .5) node[left] {$y_5 = F(q_3)$} -- (8.5, .5) %
    (-.5, 0) node[left] {$y_6 = F(p_1)$} %
    (-.5, -.5) node[left] {$y_7 = F(p_4)$} -- (8.5, -.5) %
    (-.5, -1.25) node[left] {$y_8 = F(p_3)$} -- (8.5, -1.25) ; %

    \draw[dashed] %
    (2.25, 5) -- (3.25, 5) %node[right] {$C_8$} %
    (0, 3.75) -- (3.25, 3.75) %node[right] {$C_7$} %
    (0, 2) -- (6.5, 2) %node[right] {$C_6$} %
    (0, 1.25) -- (6.5, 1.25) %node[right] {$C_5$}       %
    (0, .75) -- (8, .75) %node[right] {$C_4$} %
    (0, .25) -- (8, .25) %node[right] {$C_3$} %
    (4.5, -.25) -- (8, -.25) %node[right] {$C_2$} %
    (4.5, -.875) -- (5.75, -.875) ; %node[right] {$C_1$} ;

    \draw[red] (0, 0) node {$\bullet$} -- (0, 4.75) %
    (2.25, 1.5) node {$\bullet$} -- (2.25, 5.25) %
    (5.75, -1.25) node {$\bullet$} -- (5.75, 2.5) %
    (8, -.5) node {$\bullet$} -- (8, 1) ; %

    \draw[blue] (1.5, 4.75) node {$\times$} -- (1.5, 1.5) %
    (3.25, 5.25) node {$\times$} -- (3.25, .5) %
    (4.5, .5) node {$\times$} -- (4.5, -1.25) %
    (6.5, 2.5) node {$\times$} -- (6.5, -.5) ;
    
    \filldraw[white,draw=black] (0, 1.25) circle (1.5pt);
    \draw[decorate,decoration={brace,raise=2pt}]
    (0, 2.5) -- (0, 4.75); %
    \draw (-.25, 3.625) node [left] {$\mu^{(p)}_{1, 2}$}; %
    \filldraw[white,draw=black] (3.25, 1.25) circle (1.5pt);
    \draw[decorate,decoration={brace,raise=2pt}] (6.5, 2.5) -- (6.5, 1.5); %
    \draw (6.75, 2) node [right] {$\mu^{(q)}_{4, 3}$};
    
  \end{tikzpicture}
  \caption[Example of construction of chains for a problem with
    general real masses]%
  {Example of construction of chains for a problem with general masses
    (color online).  Red points and blue crosses mark the values of
    the cumulative distribution function~$F$ at supply points~$p_i$
    and demand points~$q_j$ according to the convention of left
    continuity. %

    Small white circles represent a pair of neighboring demand and
    supply elements. Chains connecting some neighboring mass elements
    are shown with dashed lines. %
    All chains have the same structure in each horizontal stratum
    delimited with dotted lines.  Capacity $m_k := y_{k - 1} - y_k$ of
    stratum~$k$ measures the amount of mass exchanged in that stratum.
    For example, the subsegment denoted $\mu^{(p)}_{1, 2}$ (resp.~$\mu^{(q)}_{4, 3}$)
    represents the share of supply located at~$p_1$ (resp.\ of demand located at~$q_4$) 
    that participates in the mass exchange in stratum~$2$ (resp.~$3$).  Detailed explanations are given
    in the text. %

    Observe that the problem is unbalanced, and chains in strata $5$
    and~$6$ have three supplies and two demands.}
  \label{fig:chainexnonint}
\end{figure}

Infinitesimal elements of supply and demand are pairs of the form
$(p_i, y')$ with $F(p_i) \le y'\le F(p_i + 0)$ and $(q_j, y'')$ with
$F(q_i + 0) \le y''\le F(q_i)$.  Geometrically a supply element $(p_i,
y')$ (demand element $(q_j, y'')$) corresponds to the point $(p_i,
y')$ (respectively, $(q_j, y'')$) in the vertical segment
corresponding to the supply~$p_i$ (demand~$q_j$) in the graph of the
cumulative distribution function~$F$ (see
Figure~\ref{fig:chainexnonint}).

For an infinitesimal element of supply $(p_i, y)$ define
\begin{gather*}
  r(p_i, y) = \min\{q_j \in Q\colon q_j > p_i,
  F(q_j + 0) \le y\le F(q_j)\},\\
  \ell(p_i, y) = \max\{q_j \in Q\colon q_j < p_i,
  F(q_j + 0) \le y\le F(q_j)\}
\end{gather*}
(with the usual convention $\min \varnothing = \infty$, $\max
\varnothing = -\infty$) and call the mass elements $(r(p_i, y), y)$
and~$(\ell(p_i, y), y)$ respectively the \emph{right neighbor} and the
\emph{left neighbor} of~$(p_i, y)$ if $r(p_i, y)$ and~$\ell(p_i, y)$
are finite.  The definition of right and left neighbors is then
extended to elements of demand by defining $r(q_j, y) = p_i$ whenever
$q_j = \ell(p_i, y) > -\infty$ and $\ell(q_j, y) = p_i$ whenever $q_j
= r(p_i, y) < \infty$.  Inspection of Figure~\ref{fig:chainexnonint}
should make these definitions clear.

\newcounter{tempii}
\setcounter{tempii}{\value{theorem}}
\setcounter{theorem}{\value{tempi}}
\begin{definition}[real-valued case]
  \label{def:chain-real}
  A \emph{chain} is a sequence of elements of mass that has one of the
  forms
  \begin{enumerate}
  \item $((p_{i_1}, y), (q_{j_1}, y), \dots, (p_{i_k}, y), (q_{j_k}, y))$ 
    with $\ell(p_{i_1}, y) = -\infty$, $r(q_{j_k}, y) = \infty$;
  \item $((q_{j_0}, y), (p_{i_1}, y), \dots, (q_{j_{k - 1}}, y), (p_{i_k}, y))$
    with $\ell(q_{j_0}, y) = -\infty$, $r(p_{i_k}, y) = \infty$;
  \item $((p_{i_1}, y), (q_{j_1}, y), \dots, (q_{j_{k - 1}}, y),
    (p_{i_k}, y))$ with $\ell(p_{i_1}, y) = -\infty$, $r(p_{i_k}, y) =
    \infty$.
  \end{enumerate}
  Here $k\ge 1$ and $q_{j_{m - 1}} = \ell(p_{i_m}, y) > -\infty$,
  $q_{j_m} = r(p_{i_m}, y) < \infty$ for all $m$ between $1$ and~$k$
  except the cases specified above.
\end{definition}
\setcounter{theorem}{\value{tempii}}

Note that chains have similar structure inside \emph{strata} defined
in the above graphical representation as bands separated by horizontal
lines corresponding to ordinates from the set $\{F(p_1 \pm 0),
\dots,\allowbreak F(p_M \pm 0),\allowbreak F(q_1 \pm 0), \dots, F(q_N
\pm 0)\}$: within each stratum all left and right neighbors are the
same and only the $y$ parameters differ.

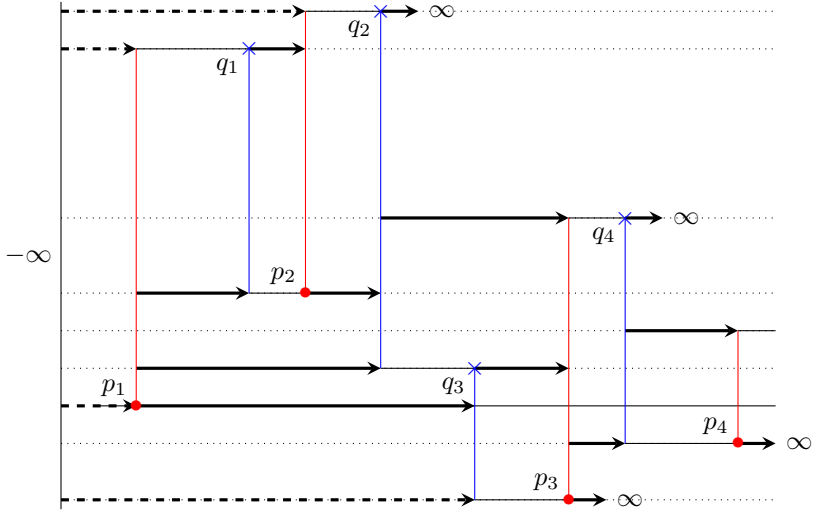
\begin{figure}
  \centering
  \begin{tikzpicture}
    \draw (-.5, 0) -- (8.5, 0) ;

    \begin{scope}[->,>=stealth,very thick]
      \draw (0, 0) -- (4.5, 0) ; %
      \draw (0, .5) -- (3.25, .5) ; %
      \draw (0, 1.5) -- (1.5, 1.5) ; %
      \draw (2.25, 1.5) -- (3.25, 1.5) ; %
      \draw (5.75, -.5) -- (6.5, -.5) ; %
      \draw (5.75, -1.25) -- (6.25, -1.25) node [right] {$\infty$} ; %
      \draw (8, -.5) -- (8.5, -.5) node [right] {$\infty$} ; %
      \draw (1.5, 4.75) -- (2.25, 4.75) ; %
      \draw (3.25, 5.25) -- (3.75, 5.25) node [right] {$\infty$} ; %
      \draw (3.25, 2.5) -- (5.75, 2.5) ; %
      \draw (4.5, .5) -- (5.75, .5) ; %
      \draw (6.5, 2.5) -- (7, 2.5) node [right] {$\infty$} ; %
      \draw (6.5, 1) -- (8, 1) ; %
    \end{scope}

    \draw[dotted] %
    (-1, 5.25) -- (8.5, 5.25) %
    (-1, 4.75) -- (8.5, 4.75) %
    (-1, 2.5) -- (8.5, 2.5) %
    (-1, 1.5) -- (8.5, 1.5) %
    (-1, 1) -- (8.5, 1) %
    (-1, .5) -- (8.5, .5) %
    (-1, -.5) -- (8.5, -.5) %
    (-1, -1.25) -- (8.5, -1.25) ; %

    \draw (-1, -1.375) -- (-1, 5.375) ;
    \draw (-1, 2) node [left] {$-\infty$} ;
    \begin{scope}[->,>=stealth,dashed,very thick]
      \draw (-1, 5.25) -- (2.25, 5.25) ; %
      \draw (-1, 4.75) -- (0, 4.75) ; %
      \draw (-1, 0) -- (0, 0) ; %
      \draw (-1, -1.25) -- (4.5, -1.25) ; %
    \end{scope}

    \draw (0, 4.75) -- (1.5, 4.75) (1.5, 1.5) -- (2.25, 1.5)
    (2.25, 5.25) -- (3.25, 5.25) (3.25, .5) -- (4.5, .5)
    (4.5, -1.25) -- (5.75, -1.25) (5.75, 2.5) -- (6.5, 2.5)
    (6.5, -.5) -- (8, -.5) (8, 1) -- (8.5, 1) ;

    \draw[red]
    (0, 0) node {$\bullet$} node [black, above left] {$p_1$} 
    -- (0, 4.75) %
    (2.25, 1.5) node {$\bullet$} node [black, above left] {$p_2$}
    -- (2.25, 5.25) %
    (5.75, -1.25) node {$\bullet$} node [black, above left] {$p_3$}
    -- (5.75, 2.5) %
    (8, -.5) node {$\bullet$} node [black, above left] {$p_4$}
    -- (8, 1) ; %

    \draw[blue]
    (1.5, 4.75) node {$\times$} node [black, below left] {$q_1$}
    -- (1.5, 1.5) %
    (3.25, 5.25) node {$\times$} node [black, below left] {$q_2$}
    -- (3.25, .5) %
    (4.5, .5) node {$\times$} node [black, below left] {$q_3$}
    -- (4.5, -1.25) %
    (6.5, 2.5) node {$\times$} node [black, below left] {$q_4$}
    -- (6.5, -.5) ;
    
  \end{tikzpicture}
  \caption{Lists $\mathcal{L}$ (solid arrows) and $\mathcal{L}_0$
    (dashed arrows) encoding the structure of the
    histogram from Figure~\protect\ref{fig:chainexnonint}. See
    explanations in the text.}
  \label{fig:listex}
\end{figure}

\subsection{Data structure and algorithm for computing chains}
\label{sec:data-struct-algor}

We now describe how to efficiently compute and store the structure of
chains and strata for a given histogram.  This discussion applies both
for the real and unitary case (the latter is degenerate in that all
elements of each supply and demand point belong to a single stratum,
\textit{cf} Figures~\ref{exchain} and~\ref{fig:chainexnonint}).  Our
construction is an adaptation of that of Aggarwal et al
\cite[Section~3]{aggarwal1992efficient} with somewhat different
terminology and notation.

The basic storage structure can be described as follows.  Observe that
for a supply point~$p_i$ the function $r(p_i, \cdot)$ is piecewise
constant and right continuous on the segment $[F(p_i), F(p_i + 0)]$.
For each $p_i$ build a list consisting of triples $(p_i, y'_{i,m},
r(p_i, y'_{i.m}))$ in the increasing order of~$m \ge 0$, where
$y'_{i,0} = F(p_i)$ and $y'_{i.m}$ corresponds to $m$th jump
of~$r(p_i, \cdot)$ as the second argument increases.  For a demand
point~$(q_j, d_j)$ build a similar list of triples $(q_j, y''_{j,m},
r(q_j, y''_{j,m}))$ where $y''_{j,0} = F(q_j)$ and $y''_{j,m}$
decreases with~$m$.  Finally build a list~$\mathcal{L}$ as
concatenation of these lists for all supply and demand points in
$P\cup Q$ in the increasing order of the abscissa.  In
Figure~\ref{fig:listex}, which features the same histogram as
Figure~\ref{fig:chainexnonint}, the elements of the combined
list~$\mathcal{L}$ are represented with thick solid arrows.  Their
order corresponds to traversing the $p_i$'s and~$q_j$'s left to right
and for each of these points, to listing the right neighbors in the
increasing order of~$y$ for $p_i$ and in the decreasing order of~$y$
for~$q_j$: in short, to traversing the continuous broken line formed
by the graph of~$F$ together with the red and blue vertical segments.

Note that all elements in~$\mathcal{L}$ that start with~$p_i$ have one
of the two following forms: $(p_i, F(p_i), q_j)$ with $q_j = r(p_i,
F(p_i))$ or $(p_i, F(q_j + 0), q_j)$ with $p_i = \ell(q_j, F(q_j +
0))$.  Similarly, elements starting with~$q_j$ have either the form
$(q_j, F(q_j), p_i)$ with $p_i = r(q_j, F(q_j))$ or $(q_j, F(p_i + 0),
p_i)$ with $q_j = \ell(p_i, F(p_i + 0))$.  Therefore all elements
of~$\mathcal{L}$ involve one of the values $F(p_i \pm 0)$ or~$F(q_j
\pm 0)$ and hence $\mathcal{L}$ has at most $2(M + N)$ elements.  To
see this refer to Figure~\ref{fig:listex} and observe, e.g., that the
function $r(p_i, \cdot)$ has a jump at~$y$ only when, during the
upward scan of the vertical segment corresponding to supply~$p_i$, one
encounters on the right the bottom end of a vertical segment
corresponding to $q_j = r(p_i, y)$ (i.e., the point with~$y = F(q_j +
0)$).  A similar observation holds for downward scan of segments
corresponding to demand elements.

The list~$\mathcal{L}$ can be regarded as a ``dictionary'' that allows
to look up the right neighbor of any supply element~$(p, y)$ or demand
element~$(q, y)$.  To do this, e.g., for $(p, y)$, locate
in~$\mathcal{L}$ an element $(\bar p, \bar y)$ immediately
preceding~$(p, y)$ and return the element $(r(\bar p, \bar y),
y)$. Again, inspection of Figure~\ref{fig:listex} should convince the
reader that this procedure is correct. Note that the search operation
in an ordered list of length~$O(M + N)$ requires an $O(\log(M + N))$
number of comparisons.

The list $\mathcal{L}$ can be built in a linear number of operations
$O(M + N)$ using the following algorithm.  Here $\mathcal{S}_p$,
$\mathcal{S}_q$ are stacks storing pairs of the form $(r, X)$ where $r
\in P\cup Q$ and $X \in \rset$.

\medskip
\begin{algorithm}\label{alg:1}
\begin{itemize}
\item Set $\mathcal{S}_p\leftarrow \varnothing$, $\mathcal{S}_q
  \leftarrow \varnothing$, list $\mathcal{L}\leftarrow \varnothing$,
  $f \leftarrow S - D$, $p \leftarrow \max P$, $q \leftarrow \max Q$;
\item loop A:
  \begin{itemize}
  \item if $p = -\infty$ and $q = -\infty$ then break loop~A;
  \item else if $p > q$ then
    \begin{itemize}
    \item set $s \leftarrow \text{supply value of~$p$}$, $P \leftarrow
      P \setminus \{p\}$;
    \item loop B: 
      \begin{itemize}
      \item if $\mathcal{S}_q = \varnothing$ then prepend $(p, f - s,
        \infty)$ to~$\mathcal{L}$ and break loop~B;
      \item pop the pair $(q', f')$ from stack~$\mathcal{S}_q$;
      \item if $f' \le f - s$ then prepend $(p, f - s, q')$
        to~$\mathcal{L}$, push the pair $(q', f')$ on stack
        $\mathcal{S}_q$ if $f' < f - s$, and break loop~B;
      \item else prepend $(p, f', q')$ to~$\mathcal{L}$;
      \end{itemize}
    \item repeat loop B;
    \item push the pair $(p, f)$ on stack~$\mathcal{S}_p$ and set $f
      \leftarrow f - s$, $p \leftarrow \max P$;
    \end{itemize}
  \item else if $p < q$ then
    \begin{itemize}
    \item set $d \leftarrow \text{demand value of~$q$}$, $Q \leftarrow
      Q \setminus \{q\}$;
    \item loop C:
      \begin{itemize}
      \item if $\mathcal{S}_p = \varnothing$ then prepend $(q, f + d,
        \infty)$ to~$\mathcal{L}$ and break loop~C;
      \item pop the pair $(p', f')$ from stack~$\mathcal{S}_p$;
      \item if $f' \ge f + d$ then prepend $(q, f + d, p')$
        to~$\mathcal{L}$, push $(p', f')$ on stack~$\mathcal{S}_p$ if
        $f' > f + d$, and break loop~C;
      \item else prepend $(q, f', p')$ to~$\mathcal{L}$;
      \end{itemize}
    \item repeat loop C;
    \item push the pair $(q, f)$ on stack~$\mathcal{S}_q$ and set $f
      \leftarrow f + d$, $q \leftarrow \max Q$;
    \end{itemize}
  \item end if;
  \end{itemize}
\item repeat loop A;
\item stop.
\end{itemize}
\end{algorithm}
\medskip

\noindent Observe that if $f$ is initialized with $S - D = F(\infty)$,
then at the exit of loop~A it will contain $F(-\infty) = 0$.  However
it is possible to initialize~$f$ with any other value, e.g.~$0$, in
which case its exit value will be smaller exactly by the amount~$S -
D$.  It is therefore not necessary to compute this quantity
beforehand.

To find leftmost mass elements of chains we also need a list
$\mathcal{L}_0$ of a similar format that stores ``right neighbors
of~$-\infty$.''  To build this list, a variant of the above procedure
is used.  While the list~$\mathcal{L}$ was built by ``prepending''
elements, i.e., adding them in front of the list, the following
algorithm uses both prepending and appending, i.e.\ adding new
elements at the end of the list.  The stacks $\mathcal{S}_p$,
$\mathcal{S}_q$ and the variable $f$ are assumed to be in the same
state as at the end of loop~A, in particular the stacks contain
exactly those $p$ and~$q$ points whose corresponding vertical segments
are ``visible from~$-\infty$.''

\medskip
\begin{algorithm}\label{alg:2}
\begin{itemize}
\item Set lists $\mathcal{L}_0 \leftarrow \varnothing$, $\mathcal{L}'
  \leftarrow \varnothing$, $\mathcal{L}'' \leftarrow \varnothing$;
\item repeat until $\mathcal{S}_q \ne \varnothing$:
  \begin{itemize}
  \item pop the pair $(q', f')$ from stack~$\mathcal{S}_q$ and append
    $(-\infty, f', q')$ to~$\mathcal{L}''$;
  \end{itemize}
\item repeat until $\mathcal{S}_p \ne \varnothing$:
  \begin{itemize}
  \item pop the pair $(p', f')$ from stack~$\mathcal{S}_p$ and prepend
    $(-\infty, f', p')$ to~$\mathcal{L}'$;
  \end{itemize}
\item if $\mathcal{L}' = \varnothing$ then
  \begin{itemize}
  \item set $(-\infty, q', f') \leftarrow$ the first element
    of~$\mathcal{L}''$ and append $(-\infty, f, q')$
    to~$\mathcal{L}_0$;
  \end{itemize}
\item else
  \begin{itemize}
  \item set $(-\infty, p', f') \leftarrow$ the last element
    of~$\mathcal{L}'$;
  \item if $\mathcal{L}'' = \varnothing$ then append $(-\infty, f,
    p')$ to~$\mathcal{L}_0$;
  \item else
    \begin{itemize}
    \item set $(-\infty, q', f') \leftarrow$ the first element
      of~$\mathcal{L}''$;
    \item append $(-\infty, f, \min\{p', q'\})$ to~$\mathcal{L}_0$;
    \end{itemize}
  \item end if;
  \end{itemize}
\item end if;
\item set $\mathcal{L}_0 \leftarrow$ concatenation of $\mathcal{L}'$,
  $\mathcal{L}_0$ and~$\mathcal{L}''$. 
\end{itemize}
\end{algorithm}

\medskip

Finally the list~$\mathcal{L}$ is scanned and the values $F(p_i \pm
0)$, $F(q_j \pm 0)$, which appear as second elements of its
constituent triples and define locations of the dotted lines
separating strata, are sorted in decreasing order to give the sequence
\begin{displaymath}
  y_0 > y_1 > \dots > y_K,
\end{displaymath}
where~$K$ is the number of strata, $k$-th stratum by definition lies
between $y_{k - 1}$ and~$y_k$, and $1 \le K \le M + N$ ($K = M + N =
8$ in the example of Figures~\ref{fig:chainexnonint},~\ref{fig:listex}).
This is the only stage in the process of building the data structure
that requires a superlinear number of operations, namely $O((M + N)
\log (M + N))$.  

\subsection{Chain decomposition of transport optimization}
\label{sec:chain-decomposition}

Observe that the initial transport optimization problem can be
replaced with a problem of transporting the Lebesgue measure supported
on ``red'' vertical segments (representing supply) to the Lebesgue
measure supported on their ``blue'' counterparts (representing
demand).  The cost function~$\bar c$ in the new problem is defined for
all points of these vertical segments, i.e., mass elements, but
depends only on their horizontal coordinates: $\bar c(p_i, y', q_j,
y'') = c(p_i, q_j)$.

% It is easy to see that any transport plan in the new problem that
% satisfies some obvious feasibility conditions can be projected to a
% feasible transport plan for the initial problem while preserving the
% total cost.  We now show that conversely, each feasible transport plan
% in the initial problem can be ``lifted'' to give some transport on the
% set of vertical segments.

Define the \emph{capacity} of $k$-th stratum as $m_k = y_{k - 1} -
y_k$ and the \emph{share of supply~$p_i$ (demand~$q_j$) in
  stratum~$k$} as $\mu^{(p)}_{i, k} = m_k$ if $F(p_i) \le y_k < y_{k -
  1} \le F(p_i + 0)$ (respectively, $\mu^{(q)}_{j, k} = m_k$ if
$F(q_i) \ge y_{k - 1} > y_k \ge F(q_i + 0)$) and~$0$ otherwise.  For
the vertical segments representing supply and demand graphically,
shares are equal to the lengths of their pieces contained between the
dotted lines (Figure~\ref{fig:chainexnonint}); we will use notation
$\mu^{(p)}_{i, k}$, $\mu^{(q)}_{j, k}$ for these subsegments as well.
Note that ${\sum_k \mu^{(p)}_{i, k} = s_i}$ (respectively, $\sum_k
\mu^{(q)}_{j, k} = d_j$).

\begin{definition}
  For a given histogram with supplies $(p_i, s_i)$ and demands~$(q_j,
  d_j)$ define a \emph{stratified transport plan} as the set of
  nonnegative values $(\gamma_{i, k; j, \ell})$, where $1\le i\le M$,
  $1\le j\le N$, and $1\le k, \ell\le K$, such that the following
  conditions are satisfied:
  \begin{equation}
    \label{eq:11}
    \sum_{i, k} \gamma_{i, k; j, \ell} = \mu^{(q)}_{j, \ell}\
    \text{for all $j$,~$\ell$},\quad
    \sum_{j, \ell} \gamma_{i, k; j, \ell} \le \mu^{(p)}_{i, k}\
    \text{for all $i$,~$k$}.
  \end{equation}
\end{definition}

Note that the numbers
\begin{equation}
  \label{eq:2}
  \gamma_{ij} = \sum_{k, \ell} \gamma_{i, k; j, \ell}
\end{equation}
form an admissible transport plan (i.e., all conditions~\eqref{eq:transport-conditions}
are satisfied).  We will call this plan the \emph{projection} of the
stratified plan in question.  The cost of a stratified transport plan
is defined as $\sum_{i, k, j, \ell} c(p_i, q_j)\, \gamma_{i, k; j,
  \ell}$; of course it coincides with the cost of its projection.

Conversely, let $\gamma = (\gamma_{ij})$, $1\le i\le M$, $1\le j\le N$
be an admissible transport plan; we call a stratified transport plan
that satisfies~\eqref{eq:2} a \emph{stratification} of~$\gamma$.  Any
admissible transport plan admits a non-empty set of stratifications.
Indeed, it is easy to check that e.g.\ for $\gamma_{i, k; j, \ell} =
\gamma_{ij} \mu^{(p)}_{i, k} \mu^{(q)}_{j, \ell}/s_id_j$ all
conditions~\eqref{eq:11}--~\eqref{eq:2} are satisfied.

We now prove that any optimal
transport plan in the initial problem can be ``lifted'' to a bundle of
disjoint transport plans operating in individual strata.
% to split the original transport optimization problem into 
% subproblems where mass is only
% moved inside each stratum~$k$.
Therefore to solve the transport optimization problem for histograms
with general real values of supply and demand, it suffices to split the
problem into transportation problems inside strata, where they reduce
to the unitary case because the mass exchanged in each stratum equals
its capacity, and solve these problems one by one.  

\begin{lemma}
  An \emph{optimal} transport plan~$\bar\gamma$ admits a
  stratification $(\bar\gamma_{i, k; j, \ell})$ that satisfies
  $\bar\gamma_{i, k; j, \ell} = 0$ whenever $\ell \neq k$
\end{lemma}

%\begin{proof}
\noindent{\em Proof:}
Indeed, let $(\gamma_{i, k; j, \ell})$
be any stratification of~$\bar\gamma$ and suppose that $\gamma_{i_0,
  k_0; j_0, \ell_0} > 0$ with $\ell_0 \ne k_0$.  Without loss of
generality we restrict the argument to the case $p_{i_0} < q_{j_0}$.

Suppose first that $\ell_0 > k_0$, i.e., that the demand
subsegment~$\mu^{(q)}_{j_0, \ell_0}$ occupies a lower stratum than the
supply subsegment~$\mu^{(p)}_{i_0, k_0}$.  The total supply located
\emph{between} these subsegments, i.e., the sum of all
$\mu^{(p)}_{i_0, k}$ with $k < k_0$ and $\mu^{(p)}_{i, k}$ with
$p_{i_0} < p_i < q_{j_0}$, is then smaller than the total demand
between these subsegments, i.e., the sum of all $\mu^{(q)}_{j, \ell}$
with $p_{i_0} < q_j < q_{j_0}$ and all $\mu^{(q)}_{j_0, \ell}$ with
$\ell < \ell_0$.  (From inspection of Figure~\ref{fig:chainexnonint} it
should be easy to see that their difference is equal to $\sum_{k_0 \le
  s < \ell_0} m_s$, although we will not need this quantity here.)
Since the first condition~\eqref{eq:11} must be fulfilled for all
$j$,~$\ell$, it follows that some demand share~$\mu^{(q)}_{j', \ell'}$
located between $\mu^{(p)}_{i_0, k_0}$ and~$\mu^{(q)}_{j_0, \ell_0}$
in the just defined sense must be satisfied with supplies located outside.
But this leads to crossing of the corresponding trajectories
(\textit{cf} Lemma~\ref{lemma:non-crossing-rule}), which implies that
the total cost of the plan $(\gamma_{i, k; j, \ell})$ can be at least
preserved, or even reduced, by a suitable rescheduling of mass
elements.

Now suppose that $\ell_0 < k_0$.  This implies the existence of extra
supply~$\mu^{(p)}_{i', k'}$ between $\mu^{(p)}_{i_0, k_0}$
and~$\mu^{(q)}_{j_0, \ell_0}$.  If this supply share is matched, it
has to feed some demand located outside, which again leads to crossing
and can be ruled out just as above.  If this supply share is not
matched (which may happen in an unbalanced problem), then
a nonzero part of the demand share $\mu^{(q)}_{j_0, \ell_0}$
can be rematched to this supply share, which is associated with the
point $p_{i'}$ located closer to~$q_{j_0}$ than $p_{i_0}$,
thus reducing the total cost.
In all cases we have a contradiction with the original assumption. $\hfill \square\newline\newline$
%\end{proof}

Note that strata are defined for piecewise constant cumulative
distributions. A simple way to apply these results to continuous
distributions of supplies and demands consists in approximate them by
piecewise constant function using, e.g. quantization techniques. 
%Approximated optimal transport plan for continuous
%densities can be computed using our method using quantization technics
%that provide relevent piecewise constant approximations of functions.

%\textsc{[It remains to discuss the mutual constraints between
%  strata.]}

%\bigskip

%\hrule

%\bigskip

%%%%%%%%%% ---------- %%%%%%%%%%
%[\textsc{From previous version}]
%%%%%%%%%% ---------- %%%%%%%%%%

%--------------------------------------------------------------------
%               SECTION V
%--------------------------------------------------------------------
\section{Practical considerations}\label{sec:prac}
In this section, we present some ways to optimize the use of the local
matching indicators in Algorithm~\ref{alg:3}.
\subsection{Exposed points}
%We state it since it permits to detect
%exposed points, hence, it can be used to save computational time.
Before applying Algorithm~\ref{alg:3}, one can detect possible
unmatched points using the following result.

\begin{lemma}[``isolation rule'']
\label{lemma:isolationrule}
Suppose that $g$ is strictly monotone and that a point  $p_i$ of  the unbalanced
chain~\eqref{eq:alter} is unmatched in
an optimal transport plan. Then if $i>1$
\begin{equation*}%\label{unmatched1}
c(p_i,q_{i-1}) \ge c(p_{i-1},q_{i-1}) 
\end{equation*}
and if $i<N$
\begin{equation*}%\label{unmatched2}
c(p_i,q_{i}) \ge c(p_{i+1},q_{i}).
\end{equation*}
 \end{lemma}

\noindent{\em Proof:} Suppose that $\sigma$ is optimal and assume for instance that $i>1$ and $c(p_i,q_{i-1}) <
   c(p_{i-1},q_{i-1})$. Thanks to Lemma~\ref{lemma:unma-isol},
 $p_i$ is not exposed, and consequently 
   $\sigma^{-1}(i-1)\leq i-1$. Thus,  $c(p_i,q_{i-1}) <
   c(p_{i-1},q_{i-1})\le c(p_{\sigma^{-1}(i-1)},q_{i-1})$. It is then
   cheaper to exclude  $p_{\sigma^{-1}(i-1)}$ and match $p_i$ with
   $q_{i-1}$, which contradicts the optimality of $\sigma$. $\hfill \square\newline\newline$
%\end{proofof}

\subsection{About the implementation and the complexity}
The cost of the algorithm can be estimated through the number of
additions and evaluations of the cost function that are required to
terminate the algorithm. These operations are only carried out in
Step~\ref{costlystep}, when computing the indicators. This section
aims at giving details about efficient ways to implement this step and
about the  complexity of the resulting procedure. 
%We focus on the most expensive part of the algorithm, namely
%Step~\ref{costlystep}, whose computational cost can be 
%reduced  for practical purposes by introducing a table that save
%and update the indicators' values that have already been
%computed.
\subsubsection{Implementation through a table of indicators}\label{sec:table}
In this section, we define a table that collects the values of indicators and then describe a
way to update it, when a negative indicator has been found. The aim of
this structure is to avoid redundant computations. We present it in the balanced case
(see~\eqref{balanced}). 

Consider a table of $N-1$ lines, where the $k$-th line
corresponds to the values
of the indicators of order $k$:
$I^p_k(1),I^q_k(1),\dots,I^q_k(N-k-1),I^p_k(N-k)$. 
%The line $k$ has
%$2(N-k)-1$ entries corresponding to the $N-k$ values of the indicators
%$I^p_k$ and the $N-k-1$ values of the indicators $I^q_k$. 
At the
beginning of the algorithm, the 
table is empty and Step~\ref{costlystep} consists in filling the line
$k$ of the table. Let us explain how to modify the table in case a
negative indicator has been found.

Following the assumptions of
Theorem~\ref{ThLMI}, consider the case where all the indicators that have been computed
currently are positive except the last one. Suppose that this one is of the
form $I^p_{k_0}(i_0)$. According to Step~\ref{removingstep}, $k_0$ pairs of supply and demand have to be matched
and removed from the current list of points $\mathcal{P}$.  Note that
the indicators that only deal with points in
$\{p_1,q_1,\dots,p_{i_0-1},q_{i_0-1},p_{i_0}\}$
% (which is defined as the empty
%set for $i_0\leq 1$) 
or in
$\{q_{i_0+k_0},p_{i_0+k_0+1},q_{i_0+k_0+1},\dots,p_N,q_N\}$ 
%(which is defined as the empty
%set for $i_0+k_0\geq N-1$) 
are not affected by this withdrawal, except that they may be
renamed. Consequently, at an order $k\leq k_0$,
$\max(0,2(i_0-k-1))+\max(0,2(N-i_0-k_0-k))$ indicators' values are
already 
known and in the line $k$ of the new table,
$2(N-k)+1-\max(0;2(i_0-k-1))-\max(0;2(N-i_0-k_0-k))$ values remain to be
computed. 
In the case the first negative indicator is of the form
$I^q_{k_0}(i_0)$, a similar reasoning shows
that in the line $k$ of the new table,
$2(N-k)+1-\max(0;2(i_0-k+1)-1)-\max(0;2(N-i_0-k_0-k)-1)$ values remain to compute. 

%On the other hand, the new
%table associated to the updated $\mathcal{P}$ has $k_0$
%lines less, and each line $k<N-k_0-1$ of
%the new table has $2(N-k)+1$ entries. 
%We deduce from this reasoning
%that in this case, only $2(N-k)+1-\max(0,2(i_0-k+1))-\max(0,2(N-i_0-k_0-k))$ 
% new values have to be computed.
%only $\min(i_0-k_0+k,2N-1-k)-\max(i_0-k_0-k,k+1)$ values
%are not valid any more since the corresponding indicator involves
%points that have been modified. Other values are not affected by the
%withdrawal. 
%Consequently, a part of
%computations in Step~\ref{costlystep} has not to be redone.
\subsubsection{Bounds for the complexity}
In the vein of the previous section, we assume up to now that all the numerical values computed during
the algorithm are saved. In this framework and as in any assignment
problem, the number of evaluations of the cost  
function cannot exceed $\frac{N(N+1)}2$. \\
The most favorable case consists in finding a negative indicator at each
step of the loop. In this case, all points are removed through
indicators of order 1. This case requires $O(N)$ additions and
evaluations of the cost function.  \\
On the opposite, the worst case corresponds to the case where all the
indicators are positive. In such a situation, no pairs are removed until the
table is full. All possible transport costs $c(p_i,q_j)$ are computed. 
Consequently, this case requires $\frac{N(N+1)}2$ evaluations of
the cost function. The number of additions is also bound by $O(N^2)$
as stated in the next theorem.
\begin{theorem}\label{th:complexity}
Denote by $C^+(N)$ the number of additions required to compute
an optimal transport plan between $N$ supplies and $N$
demands with Algorithm~\ref{alg:3}. One has: 
\begin{equation*}%\label{eq:C+}
C^+(N)\leq 3 N^2-6N.
\end{equation*}
\end{theorem}
The proof of this  result is given in Appendix.

\begin{figure}[h!]
\center
\psfrag{g}[c][t]{In-line additions}
\psfrag{N}[c][b]{Number of pairs of points $N$}
%\psfrag{data1}[l][c]{: Best case $N-1$}
\psfrag{data1}[l][c]{\tiny: $g(x)=|x|^{10^{-3}}$, $\alpha=1.22$}
\psfrag{data2}[l][c]{\tiny: $g(x)=\sqrt{|x|}$, $\alpha=1.87$}
\psfrag{data3}[l][c]{\tiny: $g(x)=|x|^{1-10^{-3}}$, $\alpha=2$}
\psfrag{data4}[l][c]{\tiny: Worst case $(N-1)^2$, $\alpha=2$}
\psfrag{100}[l][c]{\small $100$}
\psfrag{150}[c][c]{\small $150$}
\psfrag{200}[c][c]{\small $200$}
\psfrag{250}[c][c]{\small $250$}
\psfrag{300}[c][c]{\small $300$}
\psfrag{350}[c][c]{\small $350$}
\psfrag{400}[c][c]{\small $400$}
\psfrag{450}[c][c]{\small $450$}
\psfrag{500}[c][c]{\small $500$}
\psfrag{10}[c][c]{\small $10$}
%\psfrag{alpha}[l][c]{$\alpha$}
\includegraphics[width=0.8\textwidth, height=0.5\textwidth]{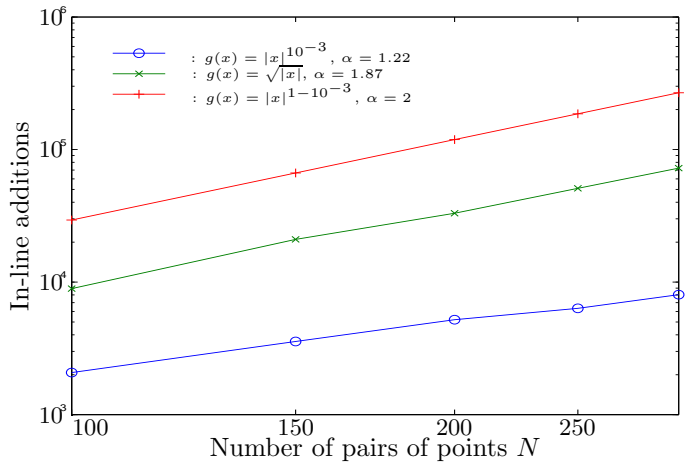}
\caption{\label{Ecp} Number of in-line additions with respect to the number of
  pairs, for various cost functions. The number $\alpha$ is the slope
  of the log-log graphs. In other words, $C^+(N) \simeq O(N^{\alpha})$.}
\end{figure}

\begin{figure}[h!]
\center
\psfrag{g}[c][t]{In-line evaluations of $g$}
\psfrag{N}[c][b]{Number of pairs of points $N$}
%\psfrag{data1}[l][c]{: Best case $N-1$}
\psfrag{data1}[l][c]{\tiny: $g(x)=|x|^{10^{-3}}$, $\alpha=1.22$}
\psfrag{data2}[l][c]{\tiny: $g(x)=\sqrt{|x|}$, $\alpha=1.88$}
\psfrag{data3}[l][c]{\tiny: $g(x)=|x|^{1-10^{-3}}$, $\alpha=2$}
\psfrag{data4}[l][c]{\tiny: Worst case $(N-1)^2$, $\alpha=2$}
\psfrag{100}[l][c]{\small $100$}
\psfrag{150}[c][c]{\small $150$}
\psfrag{200}[c][c]{\small $200$}
\psfrag{250}[c][c]{\small $250$}
\psfrag{300}[c][c]{\small $300$}
\psfrag{350}[c][c]{\small $350$}
\psfrag{400}[c][c]{\small $400$}
\psfrag{450}[c][c]{\small $450$}
\psfrag{500}[c][c]{\small $500$}
\psfrag{10}[c][c]{\small $10$}
%\psfrag{alpha}[l][c]{$\alpha$}
\includegraphics[width=0.8\textwidth, height=0.5\textwidth]{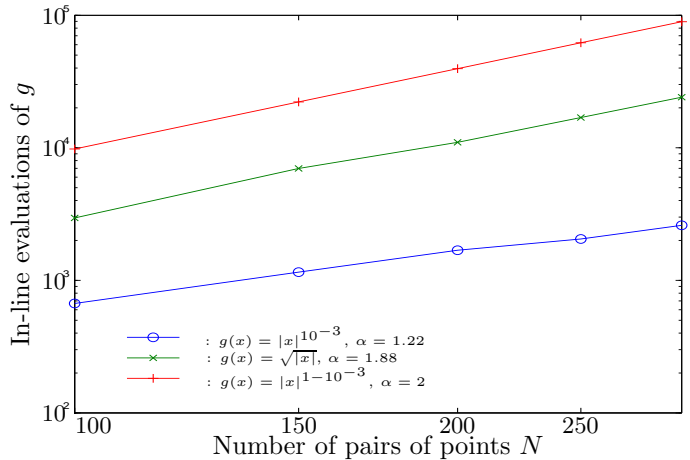}
\caption{Number of in-line evaluations of $g$ with respect to the number of
  pairs, for various cost functions. The number $\alpha$ is the slope of the log-log graphs.}\label{Ecg}
\end{figure}

In practice we observe that this upper bound is quite coarse,
especially for small values of $\alpha$. In order to better understand
this behavior, we estimated the empirical complexity of our algorithm
when $N$ is increasing, for different values of $\alpha$. 
For a fixed value of $N$, 100 samples of $N$ points are chosen
randomly in $[0,1]$, and the mean of the number of additions and
evaluations of $g$ are computed. The results are shown
in~Figures~\ref{Ecp}-\ref{Ecg} as log-log graphs. Observe that 
the less concave the cost function is, the more accurate the bound
$O(N^2)$ is. Conversely, when $\alpha$ tends towards $0$, the
complexity seems to get closer to a linear complexity. 

In order to explain this reduced complexity when $\alpha$ decreases, we
can notice that if the successive orders at which Step 4 of Algorithm~\ref{alg:3} is
visited are all bounded by $K$, then (proof in Appendix)
$$C^+(N) \leq 3N(K^2+2K+2).$$ 
Now, if we restrict ourselves to
cost functions of the type $c(x,y) = |x-y|^{\alpha}$, with $\alpha \in
(0,1]$, we also show (see Appendix) that local indicators of order $1$ tend to be
more easily negative when $\alpha$ decreases. More precisely, we show that for a given configuration of four points, the local
indicator is either positive for all $\alpha$, or there exists
$\alpha_0$, which depends only on the configuration, such that the indicator is  negative on $(0,\alpha_0]$ and
positive on $[\alpha_0,1]$.  As a consequence, the probability for an
indicator of four points to be negative increases when $\alpha$ decreases.
We conjecture that this last result remain true for all indicators (and
we checked empirically that it is). 
If this conjecture turns out to be true, the smaller $\alpha$ is, the
more probable it is for any indicator
to be negative, and the more realistic it is that the bound $K$ is
small in comparison to $N$. This would explain the previous empirical results.

\subsection{Possible improvements}%\label{sec:impro}
The use of Algorithms~\ref{alg:3}--\ref{alg:2} enables to tackle transport problems involving
real-valued histograms in $O(N^3)$ operations. Nevertheless, we emphasize that this
complexity could be reduced since there is
certainly room for improvement in the above algorithmic strategy. As
an example, identical indicators may appear in different strata and
should not be treated independently to save computational time. The
investigation of the interplay between the strata %must be investigate more accurately and
 remains for future assessment.

%{\bf HERE: preservation of hidden arcs ??}

\section*{Acknowledgments}
%This work was started during the visit of JD and JS to
%  UMI~2615 CNRS ``Laboratoire J.-V.~Poncelet,'' supported by ANR
%  through grant ANR-07-BLAN-0235 OTARIE
%  (\url{http://www.mccme.ru/~ansobol/otarie/})
This work was started during the visit of JD and JS at the Observatoire de Nice
   made possible by ANR
  through grant ANR-07-BLAN-0235 OTARIE
  (\url{http://www.mccme.ru/~ansobol/otarie/}); AS thanks the Ministry
  of National Education of France for supporting his visit to the
  Observatoire de la C{\^o}te d'Azur, where part of this text was
  written.  His work was also supported by the Russian Fund for Basic Research for the partial
  support via grant RFBR 11-01-93106-CNRSL-a and the Simons-IUM
  fellowship.
%\section*{Acknowledgments}
% Acknowledgements text here

%--------------------------------------------------------------------
%               APPENDIX
%--------------------------------------------------------------------

\section{Appendix}
\subsection*{Proof of  Theorem~\ref{th:complexity}}
Before proving Theorem~\ref{th:complexity}, let us state some intermediate
results. In what follows, we denote by $c_k^+(N)$
the number of additions required to achieve
Step~\ref{costlystep} of the algorithm for an arbitrary value of $k$.
\begin{lemma}\label{lem:c_k}
Keeping the previous notations, we have:
%\begin{equation}\label{cas:k=1}
% c_1^+(N)=3(2N-3),
%\end{equation}
%and
\begin{equation}\label{eq:c_k}
c_k^+(N)\leq 3\left(2(N-k)-1\right).
\end{equation}
%for $k>1$.
\end{lemma}
%\begin{proofof}

\noindent{\em Proof:} The proof of~\eqref{eq:c_k} in the case $k=1$ is left as exercise for the
reader. Suppose that $k>1$.
Consider for example $I^p_k(i)$ and recall that:
\begin{equation}\label{eq:idlem}
I^p_k(i)=c(p_i,q_{i+k})+\sum_{\ell=0}^{k-1}c(p_{i+\ell+1},q_{i+\ell})-\sum_{\ell=0}^{k}c(p_{i+\ell},q_{i+\ell}).
\end{equation}
The first term of this formula does not require any
addition and most of the other terms have already been computed during the
previous steps.  Indeed, the first sum
% $\sum_{\ell=0}^{k-1}c(p_{i+\ell+1},q_{i+\ell})$ 
has been computed to evaluate $I^q_{k-1}(i)$ and the second one
% sum $\sum_{\ell=0}^{k-1}c(p_{i+\ell},q_{i+\ell})$ 
has been computed to evaluate $I^p_{k-1}(i)$. It remains to
add  $c(p_{i+k},q_{i+k-1})$ to it to
compute the last sum of~\eqref{eq:idlem}. Since at given order $k$ at
most $2(N-k)-1$ indicators have to be computed, the result follows.$\hfill \square\newline\newline$
%\end{proofof}
We now consider the number of operations required between the
beginning of the algorithm and the first occurrence of
Step~\ref{removingstep}.
\begin{lemma}\label{lem:ligne}
The operations required by the algorithm between its beginning and the
first occurrence of
Step~\ref{removingstep} can be achieved with 
 $\ell^+_{k_0}(N):=3k_0(2N-k_0-2)$ additions, 
where $k_0$ denote the current value of $k$ when
Step~\ref{removingstep} occurs.
\end{lemma}
%\begin{proofof}

\noindent{\em Proof:} Between the beginning of the algorithm and the first occurrence of
Step~\ref{removingstep}, only positive indicators
have been computed, except for the current value of $k=k_0$. This means that
Step~\ref{costlystep} has been carried out for $k=1,\dots,k_0$ since the
beginning.  The corresponding number of additions is
bounded by $\sum_{k=1}^{k_0}c_k^+(N)$. Thanks to
Lemma~\ref{lem:c_k}, the result follows.$\hfill \square\newline\newline$
%\end{proofof}
Recall now that after Step~\ref{removingstep}
has been achieved, the parameter $k$ is set to $1$. The previous
arguments consequently apply to evaluate the number of additions
between two occurrences of
Step~\ref{removingstep}, i.e. between two
withdrawals. In this way, one finds that this  
number is bounded by $\ell^+_{k_0'}(N')$, where $N'$ and $k_0'$ are
the current values of $N$ and $k$ at the last occurrence of 
Step~\ref{removingstep}. Note that
$\ell^+_{k_0'}(N')$ is a coarse upper bound because we are not
considering the first occurrence of this step and a part of the
indicators has already been computed as explained in
Section~\ref{sec:table}. 

We are now in position to prove Theorem~\ref{th:complexity}.

%\begin{proofof}
\noindent{\em Proof (of Theorem~\ref{th:complexity}):} Let $k_0, k_1,\dots,k_s$
be the successive orders at which the 
Step~\ref{removingstep} of the algorithm is
visited. Observe that some of these numbers can be equal. Assume also
that only one negative indicator was found at each of these orders,
which is the worst case for complexity. As a consequence,
$\sum_{i=0}^s k_i = N$, and the number of additions required for the
whole algorithm is lower than 
\begin{equation*} 
%  \label{eq:upper_bound_addition}
  C^+ \leq \sum_{i=0}^s \ell_{k_i}^+(N-\sum_{j=0}^{i-1} k_j),
\end{equation*}
where $\ell_{k}^+$ is defined in Lemma~\ref{lem:ligne}.
%the number of additions required by the
%algorithm between its beginning, starting from $2m$ points,  and its
%$k^{th}$ occurence. 
Using Lemma~\ref{lem:ligne}, we compute
\begin{eqnarray*}
%  \label{eq:upper_bound_addition_details}
  C^+ &\leq& \sum_{i=0}^s 3k_i(2(N-\sum_{j=0}^{i-1} k_j) -k_i-2) \\
&=& \sum_{i=0}^{s-1} 3k_i(2(N-\sum_{j=0}^{i-1} k_j) -k_i -2) + 3 k_s(2(N-\sum_{j=0}^{s-1} k_j) -k_s-2)\\
&=& \sum_{i=0}^{s-1} 3k_i(2(N-\sum_{j=0}^{i-1} k_j) -k_i -2) + 3(N-\sum_{j=0}^{s-1} k_j) (N-\sum_{j=0}^{s-1} k_j - 2)\\
&=& 3N^2 - 6N - 6 \sum_{i=0}^{s-1} \sum_{j=0}^{i-1} k_ik_j - 3 \sum_{i=0}^{s-1} k_i^2 + 3 (\sum_{j=0}^{s-1} k_j)^2\\
&=& 3N^2 - 6N.
\end{eqnarray*}$\hfill \square\newline\newline$

\subsection*{Alternative complexity upper bound}
Suppose that the first occurence of Step 4 is achieved at level $k_0$,
in $3k_0(2N-k_0-2)$ additions. At this point, we remove $2k_0$ points in the total chain. Observe that the number of indicators of order $k$ that have changed after this removal of $2k_0$ points is at most $2k+1$. Let $k_1$ be the next order at which Step 4 is visited. If the indicators computed during the first pass of the algorithm have been kept in memory, this means that the number of additions necessary in the second pass is smaller than $3\sum_{k=1}^{k_1} (2k+1) = 3(k_1^2 +2 k_1)$. This yields an alternative upper bound of the whole algorithm complexity
\begin{equation}
  \label{eq:upper_bound_complexity_bis}
  C^+ \leq 3k_0(2N-k_0-2) + \sum_{i=1}^s  3(k_i^2 +2 k_i).
\end{equation}
If the successive orders at which Step 4 of the algorithm is
visited are all bounded by $K$, then $C^+ \leq 3N(K^2+K+2)$.

\subsection*{Sign of indicators for costs $|x-y|^{\alpha}$}
Consider four consecutive points $p_i$, $q_i$, $p_{i+1}$, $q_{i+1}$ in a chain. Assume without loss of generality that $|q_{i+1}-p_i| =1$. Let $a=|q_{i}-p_i|$, $b=|p_{i+1}-q_i|$ and $c=|q_{i+1}-p_{i+1}|$, so $b=1-a-c$.  
Assume that $b\leq \min(a,c)$ and define
\begin{equation}
f(\alpha) = b^{\alpha} + 1 - a^{\alpha}- c^{\alpha}. 
\end{equation}
It can be shown that if $b = 1-a-c \geq ac$, then $f$ is positive and increasing on $[0,1]$ (this result can be seen as a refined version of the rule of three). Indeed, the derivative of $f$ is
$$f'(\alpha) = \log(b)b^{\alpha}-\log(a)a^{\alpha}-\log(c)c^{\alpha}.$$
If  $b \geq ac$, then $\log(b) \geq \log(a) + \log(c)$, which implies that
$$f'(\alpha) \geq \log(a) (b^{\alpha}-a^{\alpha}) +\log(c)(b^{\alpha}-c^{\alpha})\geq 0.$$
Since $f(0)=0$, the result follows.
As a consequence, for all costs of the form $|x-y|^{\alpha}$, if $b\geq ac$, the indicator $I_1^p(i)$ will be positive.

Now, assume that $b = 1-a-c < ac$. In this case, the indicator $I_1^p(i)$ can be negative if $\alpha$ is small enough. Indeed, $f(0)=0$ and $f'(0) < 0$, which implies that $f$ is negative in the right neighborhood of $0$. Now, $f(1) = 2 -2a -2c \geq 0$, which means that the indicator is positive for $\alpha$ close to $1$. 

Consequently, there exists $\alpha_0$ such that $f(\alpha_0)=0$. Moreover, we can assume that $f'(\alpha_0)>0$, which means that $\log(b)b^{\alpha_0}-\log(a)a^{\alpha_0}-\log(c)c^{\alpha_0}>0$. Now consider $\alpha>\alpha_0$. One has successively:
\begin{eqnarray*}
f'(\alpha)&=&\log(b)b^{\alpha_0}b^{\alpha-\alpha_0}-\log(a)a^{\alpha}-\log(c)c^{\alpha}\\
&>&(\log(a)a^{\alpha_0}+\log(c)c^{\alpha_0})b^{\alpha-\alpha_0}-\log(a)a^{\alpha}-\log(c)c^{\alpha}\\
&=&(\log(a)a^{\alpha_0}+\log(c)c^{\alpha_0})b^{\alpha-\alpha_0}-\log(a)a^{\alpha}-\log(c)c^{\alpha}\\
&=&\log(a)a^{\alpha_0}(b^{\alpha-\alpha_0}-a^{\alpha-\alpha_0})+\log(c)c^{\alpha_0}(b^{\alpha-\alpha_0}-c^{\alpha-\alpha_0})\\
&>&0.
\end{eqnarray*}
This implies that if an indicator or order $1$ is negative for a given $\alpha$ in $[0,1]$, it will remain negative for smaller powers. 
%We also conjecture that this last result is true for all indicators. In this case, the smaller $\alpha$, the easier it is for any indicator to be negative.

\bibliography{transport}
\bibliographystyle{abbrv}

\end{document}